\documentclass[12pt,a4paper]{article}
\usepackage{graphicx}
\usepackage{amsfonts}
\usepackage{float}
\usepackage[latin1]{inputenc}
\usepackage{blkarray}
\usepackage{color}
\usepackage{amsmath}
\usepackage{amstext}
\usepackage{amssymb}
\usepackage{array}%,arydshln}
\usepackage{multicol}
\usepackage{pdfpages}
\usepackage{multirow}

\usepackage{threeparttable}
\usepackage{url}
\usepackage{enumerate}
\usepackage[ruled,vlined]{algorithm2e}
\usepackage{algorithmic}
\usepackage{calc}
\usepackage{caption}
\usepackage{subcaption}
\usepackage{marginnote}

\def\N{\Bbb N}
\def\R{\Bbb R}

\def\D{\cal D}

\def\<{\langle}
\def\>{\rangle}

\def\e{\varepsilon}

\def\Chi{\raise .3ex \hbox{\large $\chi$}}

\def\n{{\bf }}

\def\mass{{\rm h_2o} }
\def\energy{e}
\def\WI{W\!\!I}

\def\[{\Bigl [}
\def\]{\Bigr ]}
\def\({\Bigl (}
\def\){\Bigr )}
\def\[{\Bigl [}
\def\]{\Bigr ]}
\def\({\Bigl (}
\def\){\Bigr )}
\def\l{\iota}

\def\div{{\mbox{\rm div}}}
\def\dsp{\displaystyle}
\def\K{{\bf K}}
\def\x{{\bf x}}

\def\n{{\bf n}}
\def\s{{\bf s}}

\def\q{{\bf q}}
\def\G{{\Gamma}}

\def\D{{\cal D}}
\def\cells{{\cal M}}
\def\faces{{\cal F}}
\def\nodes{{\cal V}}
\def\edges{{\cal E}}

\def\e{{\mathfrak a}}

\def\welledge {{\mathfrak a}}

\def\l{{\ell}}
\def\g{{\rm g}}

\begin{document}

%%-----------------------------
%%      the top matter
%%-----------------------------
\title{Two-phase geothermal model with fracture network and multi-branch wells}

\author{
A. Armandine Les Landes\thanks{BRGM, 3 avenue Claude-Guillemin, BP 36009, 45060 Orl\'eans Cedex 2, France, A.ArmandineLesLandes@brgm.fr},   
D. Castanon Quiroz \thanks{Universit\'e C\^ote d'Azur, Inria, CNRS, LJAD, 
  UMR 7351 CNRS, team Coffee, Parc Valrose 06108 Nice Cedex 02, France, danielcq.mathematics@gmail.com},
L. Jeannin \thanks{STORENGY, 12 rue Raoul Nordling - Djinn - CS 70001 92274 Bois Colombes Cedex, France, laurent.jeannin@storengy.com}, 
S. Lopez\thanks{BRGM, 3 avenue Claude-Guillemin, BP 36009, 45060 Orl\'eans Cedex 2, France, s.lopez@brgm.fr}, 
R. Masson\thanks{Universit\'e C\^ote d'Azur, Inria, CNRS, LJAD, 
UMR 7351 CNRS, team Coffee, Parc Valrose 06108 Nice Cedex 02, France, roland.masson@unice.fr} 
}
\maketitle
%%-----------------------------
%%      your text
%%-----------------------------

\begin{abstract}
  This paper focuses on the numerical simulation of geothermal systems in complex geological settings. The physical model is based on two-phase Darcy flows coupling the mass conservation of the water component with the energy conservation and the liquid vapor thermodynamical equilibrium.
  The discretization exploits the flexibility of  unstructured meshes to model complex geology including conductive faults as well as complex wells.  
  The polytopal and essentially nodal Vertex Approximate Gradient scheme is used for the approximation of the Darcy and Fourier fluxes combined with a Control Volume approach for the transport of mass and energy. 
Particular attention is paid to the faults which are modelled as two-dimensional interfaces defined as collection of faces of the mesh and to the flow inside deviated or multi-branch wells defined as collection of edges of the mesh with rooted tree data structure. By using an explicit pressure drop calculation, the well model reduces to a single equation based on complementarity constraints with only one well implicit unknown. The coupled systems are solved fully implicitely at each time step using efficient nonlinear and linear solvers on parallel distributed architectures. The convergence of the discrete model is investigated numerically on a simple test case with a Cartesian geometry and a single vertical producer well. Then, the ability of our approach to deal efficiently with realistic test cases is assessed on a high energy faulted geothermal reservoir operated using a doublet of two deviated wells. 
\end{abstract}  

\section{Introduction}
\label{sec_intro}
\hspace{\parindent} 
Deep geothermal systems are often located in complex geological settings, including faults or fractures. These geological discontinuities not only control fluid flow and heat transfer, but also provide feed zones for production wells. Modeling the operation of geothermal fields and the exchange of fluids and heat in the rock mass during production requires explicitly taking into account objects of different characteristic sizes such as the reservoir itself, faults and fractures, which have a small thickness compared to the characteristic size of geological formations and wells (whose radius is of the order of a few tens of centimeters).

A common way to account for these highly constrated spatial scales is based on a reduction of dimension both for the fault/fracture and the well models. Following \cite{GRANET200135,MAE02,BMTA03,FNFM03,MJE05,RJBH06,HF08,ABH09,GSDFN,BHMS2016,NBFK2019} faults/fractures will be represented as co-dimension one manifolds coupled with the surrounding matrix domain leading to the so-called hybrid-dimensional or Discrete Fracture Matrix (DFM) models. This reduction of dimension is obtained by averaging both the equations and unknowns in the fracture width and using appropriate transmission conditions at matrix fracture interfaces. In our case, the faults/fractures will be assumed to be conductive both in terms of permeability and thermal conductivity in such a way that pressure and temperature continuity can be assumed as matrix fracture transmission conditions \cite{MAE02,BMTA03,RJBH06}. This setting has been extended to two-phase Darcy flows in \cite{BGGM14,BGJMP17} and to multi-phase compositional non-isothermal Darcy flows in \cite{Xing.ea:2017}. 

The well will be modelled as a line source defined by a 1D graph with tree structure. It will be coupled to the 3D matrix domain and to the 2D faults/fractures possibly intersecting the well using Peaceman's approach. It is a widely used approach in reservoir simulation for which the Darcy or Fourier fluxes between the reservoir and the well are discretized by a two-point flux approximation with a transmissivity accounting for the unresolved pressure or temperature singularity. This leads to the concept of well or Peaceman's index defined at the discrete level and depending on the type of cell, on the well radius and geometry and on the scheme used for the discretization. Let us refer to \cite{Peaceman78} for its introduction in the framework of a two-point cell-centered finite volume scheme on square cells, to \cite{Peaceman83} for its extension to non square cells and anisotropic permeability field and to \cite{Wolf03,Aav03,CZ09} for extensions to more general well geometries and different discretizations. The coupling with the faults/fractures is considered in \cite{Xing2018}. Let us also refer to \cite{Gjerde20} for a related approach also based on a removal of the singularity induced by the well line source but at the continuous level.  \\

This paper focuses on the liquid vapor single water component non-isothermal Darcy flow model based on mass and energy conservation equations coupled with thermodynamical equilibrium and volume balance. The extension to hybrid-dimensional models follows \cite{Xing.ea:2017} with pressure and temperature continuity at matrix fracture interfaces. The thermal well model is a simplified version of the drift flux model \cite{LIVESCU2010138,shi2005} neglecting transient terms, thermal losses and cross flow in the sense that all along the well, the well behaves either as a production or an injection well. It results that using an explicit approximation of the mixture density along the well, the well model can be reduced to a single unknown, the so-called bottom hole pressure, implicitely coupled to the reservoir.  \\

The discretization of hybrid-dimensional Darcy flow models 
has been the object of many works using cell-centered Finite Volume schemes with either
Two Point or Multi Point Flux Approximations 
\cite{KDA04,ABH09,HADEH09,TFGCH12,SBN12,AELHP152D,AELHP153D}, 
Mixed or Mixed Hybrid Finite Element methods  \cite{MAE02,MJE05,HF08},
Hybrid Mimetic Mixed Methods \cite{FFJR16,AFSVV16,GSDFN,BHMS2016}, 
and Control Volume Finite Element Methods (CVFE) \cite{BMTA03,RJBH06,MF07,HADEH09,MMB2007}. 
This article focus on the Vertex Approximate Gradient (VAG) scheme accounting for polyhedral meshes. It 
has been introduced for the discretization of multiphase Darcy flows 
in \cite{EHGM-CG-12} and extended to hybrid-dimensional models in 
\cite{BGGM14,GSDFN,tracer2016,BHMS2016,Xing.ea:2017,BHMS18,Brenner2021}.

The VAG scheme uses nodal and fracture face unknowns 
in addition to the cell unknowns which can be eliminated without any fill-in. 
Thanks to its essentially nodal nature, it leads to a sparse discretization on tetrahedral meshes which
are particulary adapted to discretize complex geological features like faults defined as a collection of faces
and slanted or multi-branch wells defined as a collection of edges with tree structure. 
Compared with other nodal approaches such as CVFE methods, the VAG scheme avoid the mixing of the control volumes at the  matrix fracture interfaces, 
which is a key feature for its coupling with a transport model. As shown 
in \cite{BGGM14} for two-phase flow problems, this allows to use a coarser mesh size at the matrix fracture interface.\\

The remainder of this paper is organized as follows. Section 2 presents the physical model describing the flow and transport in the matrix domain coupled to the fracture/fault network in the hybrid-dimensional setting. 
Section 3 presents the VAG discretization of this liquid vapor non-isothermal hybrid-dimensional model. It is based on the discrete mass and energy conservations on each control volume coupled with thermodynamical equilibrium and the sum to one of the saturations. Then, the well modelling is addressed starting with the description of the well geometry as a collection of edges defining a rooted tree data structure. The source terms connecting the well to the reservoir at each well node are based on two-point fluxes with transmissivities defined by Peaceman's indexes. The derivation of the simplified well model is detailed both for production and injection wells starting from the drift flux model. We discuss at the end of Section 3 the algorithms used to solve the nonlinear and linear systems on distributed parallel architectures at each time step of the simulation. 
Finally, to demonstrate the efficiency of our approach, we present in Section 4 two numerical tests. The first test case checks the numerical convergence of the model for a vertical production well connected to an homogeneous reservoir on a family of refined Cartesian meshes. The second test case simulates the development plan of a high enthalpy faulted geothermal reservoir with slanted production and injection wells.

\section{Hybrid-dimensional non-isothermal two-phase Discrete Fracture Model}
\label{sec_model}
This section recalls, in the particular case of a non-isothermal single-component two-phase Darcy flow model,
the hybrid-dimensional model introduced in \cite{Xing.ea:2017}. 

\subsection{Discrete Fracture Network}
\label{subsec_dfn}

Let $\Omega$ denote a bounded domain of $\R^3$  assumed to be polyhedral. 
Following \cite{MAE02,FNFM03,MJE05,GSDFN,BHMS2016} the fractures are represented as interfaces of codimension 1.  
Let $J$ be a finite set and let 
$
\overline \Gamma = \bigcup_{j\in J} \overline \Gamma_j
$  
and its interior $\Gamma = \overline \Gamma\setminus \partial\overline\Gamma$ 
denote the network of fractures $\Gamma_j\subset \Omega$, $j\in J$, such that each $\Gamma_j$ is 
a planar polygonal simply connected open domain included in a plane of $\R^3$. 
\begin{figure}[htbp!]
\begin{center}
\includegraphics[width=0.37\textwidth]{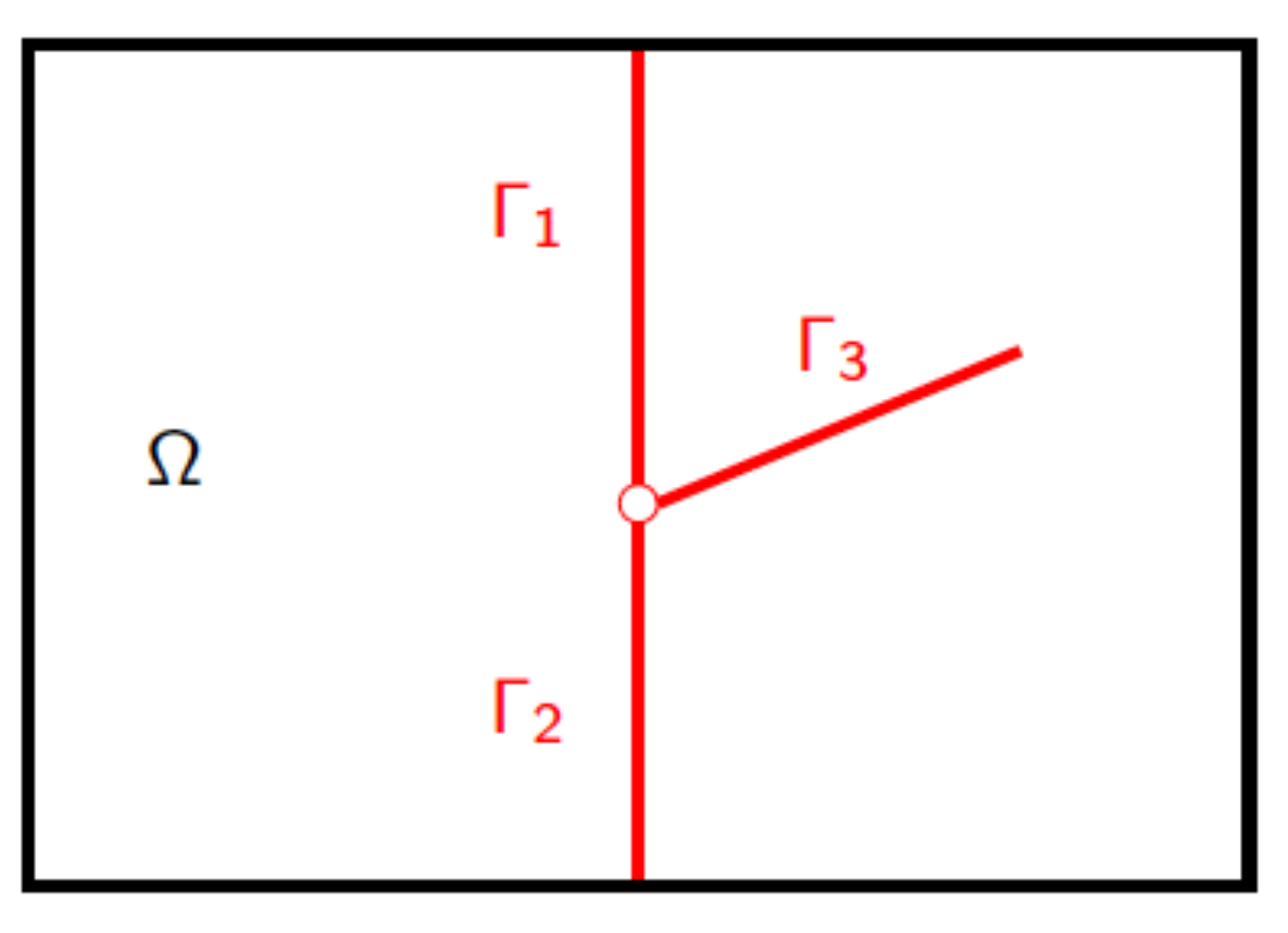}   
\caption{Example of a 2D domain with 3 intersecting fractures $\Gamma_1, \Gamma_2, \Gamma_3$.}
\label{fig_network}
\end{center}
\end{figure}
The fracture width is denoted by $d_f$ and is such that $0 < \underline d_f \leq d_f(\x) \leq \overline d_f$ 
for all $\x\in \Gamma$. 
We can define, for each fracture $j\in J$, its two 
sides $+$ and $-$. For scalar functions on $\Omega$, possibly 
discontinuous at the interface $\Gamma$
(typically in $H^1(\Omega\setminus \overline\Gamma)$), 
we denote by $\gamma^\pm$ the trace operators on the side $\pm$ of $\Gamma$. 
Continuous scalar functions $u$ at the interface $\Gamma$ (typically in $H^1(\Omega)$) 
are such that $\gamma^+ u = \gamma^- u$ 
and we denote by $\gamma$ the trace operator on $\Gamma$ for such functions. 
At almost every point of the fracture network, 
we denote by ${\bf n}^\pm$ the unit normal vector oriented outward to the side $\pm$ of $\Gamma$ 
such that ${\bf n}^+ + {\bf n}^- = 0$. For vector fields on $\Omega$,  possibly discontinuous at the interface $\Gamma$ (typically in $H_\div(\Omega\setminus\overline\Gamma)$, 
we denote by $\gamma_n^{\pm}$ the normal trace operator on the side $\pm$ of $\Gamma$ oriented w.r.t. ${\bf n}^\pm$. 

The gradient operator in the matrix domain $\Omega\setminus \overline\Gamma$ is denoted by $\nabla$ and the tangential gradient operator on the 
fracture network is denoted by  $\nabla_\tau$ such that 
$$
\nabla_\tau u  = \nabla u - (\nabla u\cdot \n^+)\n^+. 
$$
We also denote by $\div_{\tau}$ the tangential divergence operator on the fracture network, and by 
$d\tau({\bf x})$ the Lebesgue measure on $\Gamma$.

We denote by $\Sigma$ the dimension $1$ open set defined by the intersection of the fractures 
excluding the boundary of the domain $\Omega$, i.e. the interior of 
$\bigcup_{\{(j,j')\in J\times J\,|\, j\neq j'\}} \partial \Gamma_j \cap \partial \Gamma_{j'} 
\setminus \partial \Omega$.  

For the matrix domain, Dirichlet (subscript $D$) and Neumann (subscript $N$) boundary conditions are imposed on 
the two dimensional open sets $\partial \Omega_D$ and $\partial \Omega_N$ respectively where 
$\partial \Omega_D\cap \partial \Omega_N = \emptyset$, 
$\partial \Omega = \overline {\partial \Omega_D}\cup \overline {\partial \Omega_N}$.
Similarly for the fracture network, the Dirichlet and Neumann boundary conditions are imposed on the one dimensional open sets $\partial \Gamma_D$ and $\partial \Gamma_N$ respectively where 
$\partial \Gamma_D\cap \partial \Gamma_N = \emptyset$, 
$\partial \Gamma \cap \partial\Omega = \overline {\partial \Gamma_D}\cup \overline {\partial \Gamma_N}$.

%Let $\gamma_{n_{\partial \Gamma_j}}, j\in J$ denote the normal trace operator at the fracture $\Gamma_j$ boundary oriented outward to $\Gamma_j$. 

\subsection{Non-isothermal two-phase flow model}

We consider in this work a two-phase liquid gas, single water component, and non-isothermal Darcy flow model. 
The liquid ($\l$) and gas ($\g$) phases are described by their pressure $p$ (neglecting capillary effects), temperature $T$ and pore volume fractions or saturations $s^\alpha$, $\alpha\in\{\l,\g\}$.
Let us also introduce the mass fraction $c^\alpha$ of the water component in phase $\alpha$, equal to $1$ for a present phase $\alpha$ but lower than $1$ for an absent phase. It will be used below to express the thermodynamical equilibrium as complementary constraints. 

 For each phase $\alpha$, we denote by  
 $\rho^\alpha(p,T)$ its mass density, by 
 $\mu^\alpha(p,T)$  its dynamic viscosity, by 
 $e^\alpha(p,T)$ its specific internal energy, and by 
 $h^\alpha(p,T)$ its specific enthalpy. The rock energy density is denoted by $E_r(p,T)$. 

 The reduction of dimension in the fractures leading to the hybrid-dimensional model
 is obtained by integration of the conservation equations
 along  the width of the fractures complemented by transmission conditions at both sides of the
 matrix fracture interfaces (see \cite{Xing.ea:2017}). In the following, $p_m,T_m, s^\alpha_m, c^\alpha_m$ denote the pressure, temperature, saturations, and mass fractions 
 in the matrix domain $\Omega\setminus\overline\Gamma$, and $p_f,T_f, s^\alpha_f, c^\alpha_f$ are the pressure, temperature, saturations and mass fractions in the fractures
 averaged along the width of the fractures. 
 The permeability tensor is denoted by ${\bf K}_m$ in the matrix domain and we denote by ${\bf K}_f$ the tangential permeability tensor in the fractures (average value along the fracture width assuming that the permeability tensor in the fracture has the normal as principal direction).  
 The porosity (resp. thermal conductivity of the rock and fluid mixture)
 is denoted by $\phi_m$ (resp. $\lambda_m$) in the matrix domain and by $\phi_f$ (resp. $\lambda_f$) along the fracture network (average values along the fracture width).
 The relative permeability of phase $\alpha$ as a function of the phase saturation is denoted by $k^\alpha_{r,m}$ in the matrix and by $k^\alpha_{r,f}$ in the fracture network. 
 The gravity acceleration vector is denoted by ${\bf g}$. \\
 
 The set of equations couples the mass, energy and volume balance equations in the matrix 
\begin{equation}
\label{hybrid_cons_mat}
\left\{  
\begin{array}{ll}
\dsp \phi_m  ~\partial_t \(\sum_{\alpha\in \{\l,\g\}}  \rho^\alpha(p_m,T_m)s_m^\alpha c^\alpha_m\)  + \div(\q^\mass_m )  =  0,\\
\dsp \phi_m ~\partial_t \( \sum_{\alpha\in \{\l,\g\}}  \rho^\alpha(p_m,T_m)e^\alpha(p_m,T_m) s^\alpha_m c^\alpha_m\) + (1-\phi_m) \partial_t E_r(p_m,T_m)  
+\div ( \q^\energy_{m} )  = 0, \\
\dsp \sum_{\alpha\in \{\l,\g\}} s^\alpha_m = 1,
\end{array}\right. 
\end{equation}
in the fracture network 
\begin{equation}
  \label{hybrid_cons_frac}
\left\{  
\begin{array}{ll}
\dsp d_f \phi_f ~ \partial_t \(\sum_{\alpha\in \{\l,\g\}} \rho^\alpha(p_f,T_f)s^\alpha_f c^\alpha_f\)  +  \div_\tau( \q^\mass_f) - \gamma_n^+ \q^\mass_{m} - \gamma_n^- \q^\mass_{m}  = 0,\\
\dsp d_f \phi_f ~\partial_t \( \sum_{\alpha\in \{\l,\g\}}\rho^\alpha(p_f,T_f)e^\alpha(p_f,T_f) s^\alpha_f c^\alpha_f\) +  d_f (1-\phi_f) \partial_t E_r(p_f,T_f)   \\
\qquad \qquad\qquad \qquad \qquad \qquad + ~\div_\tau ( {\bf q}^\energy_{f} ) - \gamma_n^+ \q^\energy_{m} - \gamma_n^- \q^\energy_{m} = 0,\\
\dsp \sum_{\alpha\in \{\l,\g\}} s^\alpha_f = 1,
\end{array}\right. 
\end{equation}
with the thermodynamical equilibrium for $i=m,f$
\begin{equation}
  \label{thermo_eq}
\left\{\begin{array}{r@{\,\,}c@{\,\,}l}  
c^{\g}_i p_i - p_{\rm sat}(T_i) c^{\l}_i &=& 0,\\  [1ex]
\min\(s^{\l}_i, 1-c^{\l}_i \) &=& 0, \\[1ex]
\min\(s^{\g}_i, 1-c^{\g}_i \) &=& 0, 
\end{array}\right. 
\end{equation}
where $p_{\rm sat}(T)$ is the vapor saturated pressure as a function of the temperature $T$. 

The Darcy and Fourier laws provide the mass and energy fluxes in the matrix 
\begin{equation}
  \label{Darcy_Fourier_mat}
  \begin{array}{ll@{\,\,}c@{\,\,}l}  
   & \q^\mass_{m} &=& \dsp \sum_{\alpha\in \{\l,\g\}} \q_m^\alpha, \\[3ex]
   & \q_m^\alpha &=& \dsp  c^\alpha_m {\rho^\alpha(p_m,T_m)\over \mu^\alpha(p_m,T_m)} k_{r,m}^\alpha(s_m^\alpha) {\bf V}^\alpha_m, \\[2ex]
   & \q^\energy_{m} &=&   \dsp \sum_{\alpha\in \{\l,\g\}} h^\alpha (p_m,T_m) \q^\alpha_m - \lambda_m \nabla T_m,
    \end{array}
\end{equation}
and in the fracture network 
\begin{equation}
  \label{Darcy_Fourier_frac}
  \begin{array}{ll@{\,\,}c@{\,\,}l}  
    & \q^\mass_{f} &=& \dsp \sum_{\alpha\in \{\l,\g\}} \q_f^\alpha, \\[3ex]
    & \q^\alpha_{f} &=& \dsp c^\alpha_f {\rho^\alpha(p_f,T_f)\over \mu^\alpha(p_f,T_f)} k_{r,f}^\alpha(s_f^\alpha){\bf V}^\alpha_f, \\[2ex]
    & \q^\energy_{f} &=&  \dsp \sum_{\alpha\in \{\l,\g\}} h^\alpha(p_f,T_f)\q^\alpha_f - d_f \lambda_f \nabla_\tau T_f,
      \end{array}
\end{equation}
where
\begin{equation*}
{\bf V}^\alpha_m  = - {\bf K}_m \( \nabla p_m - \rho^\alpha(p_m,T_m) {\bf g}\),\quad\quad 
{\bf V}^\alpha_f  = - d_f {\bf K}_f \( \nabla_\tau p_f - \rho^\alpha(p_f,T_f) {\bf g}_\tau\),  
\end{equation*}
and ${\bf g}_\tau = {\bf g} - ({\bf g}\cdot\n^+)\n^+$.\\

The system \eqref{hybrid_cons_mat}-\eqref{hybrid_cons_frac}-\eqref{Darcy_Fourier_mat}-\eqref{Darcy_Fourier_frac} is closed with transmission conditions at the matrix fracture interface $\Gamma$.  
These conditions state the continuity of the pressure and temperature at the matrix fracture interface assuming that the fractures do not act as barrier neither for the Darcy flow nor for the thermal conductivity 
(see \cite{MAE02,FNFM03,MJE05,Xing.ea:2017}). 
\begin{equation}
\label{Hybrid_transmission}
\begin{split}
& \gamma^+ p_m= \gamma^- p_m = \gamma p_m = p_f,\\
& \gamma^+ T_m= \gamma^- T_m = \gamma T_m = T_f. 
\end{split}
\end{equation}
At fracture intersections $\Sigma$, note that we assume mass and energy flux conservation as well as the continuity of the pressure $p_f$ and temperature $T_f$.  
Homogeneous Neumann boundary conditions are applied 
for the mass $\q^\mass_{f}$ and energy $\q^\energy_{f}$ 
fluxes at the fracture tips $\partial \Gamma\setminus \partial \Omega$. 

\section{VAG Finite Volume Discretization}
\label{sec_discretization}
\subsection{Space and time discretizations}

The VAG discretization of hybrid-dimensional two-phase Darcy flows introduced in \cite{BGGM14} considers generalized polyhedral meshes of $\Omega$ in the spirit of  \cite{Eymard.Herbin.ea:2010}. Let $\cells$ be the set of cells that are disjoint open polyhedral subsets of $\Omega$ such that
$\bigcup_{K\in\cells} \overline{K} = \overline\Omega$, for all $K\in\cells$, ${\x}_K$ denotes the so-called ``center'' of the cell $K$ under the assumption that $K$ is star-shaped with respect to ${\x}_K$. 
The set of faces of the mesh is denoted by $\faces$ and $\faces_K$ is the set of faces of the cell $K\in \cells$. 
The set of edges of the mesh is denoted by $\edges$ and $\edges_\sigma$ is the set of edges of the face $\sigma\in \faces$. 
The set of vertices of the mesh is denoted by $\nodes$ and $\nodes_\sigma$ is the set of vertices of the face $\sigma$. 
For each $K\in \cells$ we define $\nodes_K = \bigcup_{\sigma\in \faces_K} \nodes_\sigma$.   

The faces are not necessarily planar. It is just assumed that for each face $\sigma\in\faces$, there exists a so-called ``center'' of the face ${{\bf x}}_\sigma \in {\sigma}\setminus \bigcup_{\e\in \edges_\sigma} \e$ such that
$
{\x}_\sigma = \sum_{\s\in \nodes_\sigma} \beta_{\sigma,\s}~\x_\s, \mbox{ with }
\sum_{\s\in \nodes_\sigma} \beta_{\sigma,\s}=1,
$
and $\beta_{\sigma,\s}\geq 0$ for all $\s\in \nodes_\sigma$; moreover
the face $\sigma$ is assumed to be defined  by the union of the triangles
$T_{\sigma,\e}$ defined by the face center ${\x}_\sigma$
and each edge $\e\in\edges_\sigma$. 
The mesh is also supposed to be conforming w.r.t. the fracture network $\G$ in the sense that for each  $j\in J$ there exists a subset $\faces_{\G_j}$ of $\faces$ such that 
$$
\overline \G_j = \bigcup_{\sigma\in\faces_{\G_j}} \overline{\sigma}.  
$$
We will denote by $\faces_\G$ the set of fracture faces 
$$
\faces_\Gamma = \bigcup_{j\in J} \faces_{\G_j}, 
$$ 
and by 
$$
\nodes_\Gamma = \bigcup_{\sigma\in \faces_\Gamma} \nodes_{\sigma}, 
$$ 
the set of fracture nodes. 
This geometrical discretization of $\Omega$ and $\G$ is denoted in the following by $\D$. 

In addition, the following notations will be used 
$$
\cells_\s = \{K\in \cells\,|\, \s\in \nodes_K\}, 
\
\cells_{\sigma} = \{K\in \cells\,|\, \sigma\in \faces_K\}, 
$$
and 
$$
\faces_{\G,\s} = \{\sigma \in \faces_\Gamma \,|\, \s\in \nodes_\sigma\}. 
$$
\vskip 0.5cm

For $N_{t_f}\in\N^*$, let us consider the time discretization 
$t^0= 0 < t^1 <\cdots < t^{n-1} < t^n \cdots < t^{N_{t_f}} = t_f$ 
of the time interval $[0,t_f]$. We denote the time steps by 
$\Delta t^n = t^{n}-t^{n-1}$ for all $n=1,\cdots,N_{t_f}$. 

\subsection{VAG fluxes and control volumes}

\begin{figure}[H]
\begin{center}
\includegraphics[width=0.4\textwidth]{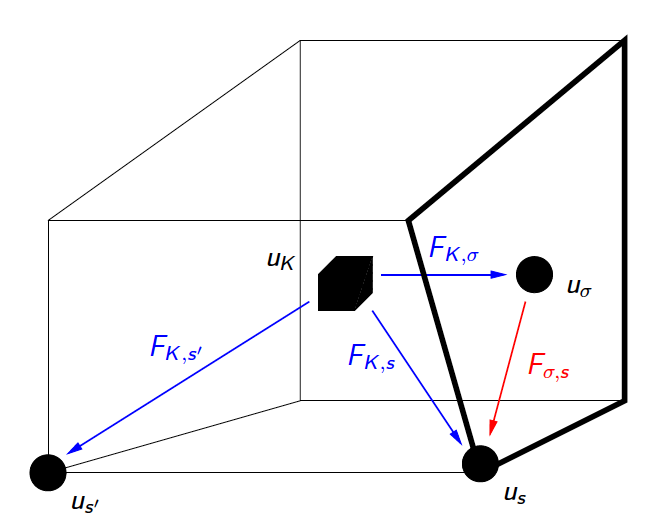}
\caption{For a cell $K$ and a fracture face $\sigma$ (in bold), 
examples of VAG degrees of freedom $u_K$, $u_\s$, $u_\sigma$, $u_{\s'}$ and 
VAG fluxes $F_{K,\sigma}$, $F_{K,\s}$, $F_{K,\s'}$, $F_{\sigma,\s}$. }
\label{fig_vag_fluxes}
\end{center}
\end{figure}

The VAG discretization is introduced in \cite{Eymard.Herbin.ea:2010} for diffusive problems 
on heterogeneous anisotropic media. Its extension to the hybrid-dimensional 
Darcy flow model is proposed in \cite{BGGM14} based upon the following vector space of degrees of freedom: 
$$
V_\D =\{v_K, v_\s, v_\sigma\in\R, K\in \cells, \s\in \nodes, \sigma\in \faces_\G\}. 
%V_\D =\{v_\nu\in\R, \nu\in \cells\cup\nodes\cup\faces_\G\}. 
$$
The degrees of freedom are illustrated in Figure \ref{fig_vag_fluxes} 
for a given cell $K$ with one fracture face $\sigma$ in bold. 

\begin{figure}[H]
\begin{center}
\includegraphics[width=0.4\textwidth]{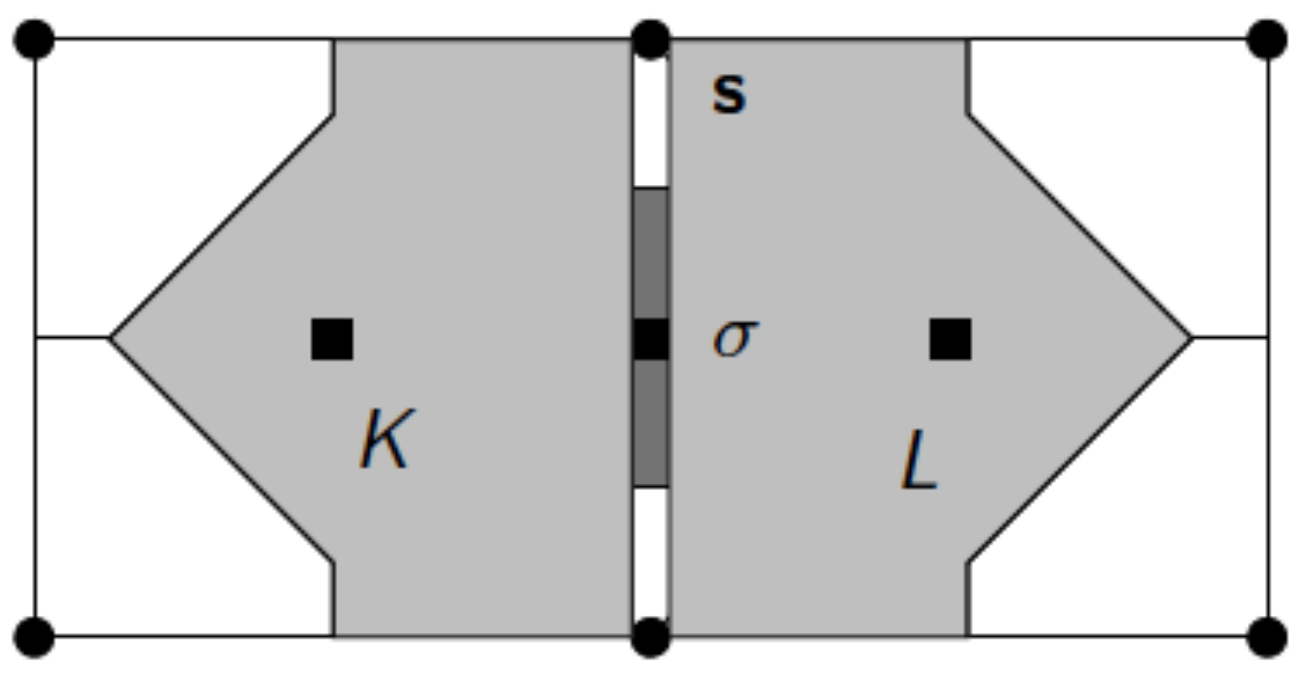}
\caption{Example of control volumes at cells, fracture face, and nodes, 
in the case of two cells $K$ and $L$ separated by 
one fracture face $\sigma$ (the width of the fracture is enlarged in this figure). 
The control volumes are chosen to avoid mixing fracture and matrix rocktypes. }
\label{fig_vag_CV}
\end{center}
\end{figure}

The matrix degrees of freedom are defined by the set of cells $\cells$ and 
by the set of nodes $\nodes\setminus\nodes_\Gamma$ excluding the nodes 
at the matrix fracture interface $\Gamma$. The fracture faces $\faces_\Gamma$ 
and the fracture nodes $\nodes_\Gamma$ are shared between the matrix and the fractures but
the control volumes associated with these degrees of freedom will belong to the fracture network (see Figure \ref{fig_vag_CV}). 
The degrees of freedom at the fracture intersection $\Sigma$ are defined by the set of nodes $\nodes_\Sigma \subset \nodes_\Gamma$ located on $\overline\Sigma$. 
The set of nodes at the Dirichlet boundaries $\overline {\partial \Omega_D}$ 
and $\overline {\partial \Gamma_D}$ is denoted by ${\cal V}_D$. 

The VAG scheme is a control volume scheme in the sense that it results, 
for each non Dirichlet degree of freedom in a mass or energy balance equation. 
The matrix diffusion tensor is assumed to be cellwise constant and the tangential 
diffusion tensor in the fracture network is assumed to be facewise constant. 
The two main ingredients are therefore the conservative fluxes and the control volumes. 
The VAG matrix and fracture fluxes are illustrated in Figure \ref{fig_vag_fluxes}. 
For $u_\D\in V_\D$, the matrix 
fluxes $F_{K,\nu}(u_\D)$ connect the cell $K\in\cells$ 
to the degrees of freedom located at the boundary of $K$, namely 
$\nu\in \Xi_K = \nodes_K \cup (\faces_K \cap \faces_\G)$. 
The fracture fluxes $F_{\sigma,\s}(u_\D)$ connect each 
fracture face $\sigma\in \faces_\G$ to its nodes $\s\in \nodes_\sigma$.  
The expression of the matrix (resp. the fracture) fluxes 
is linear and local to the cell (resp. fracture face). 
More precisely, the matrix fluxes are given by 
$$
F_{K,\nu}(u_\D) = \sum_{\nu'\in \Xi_K} T_K^{\nu,\nu'} (u_K - u_{\nu'}), 
$$
with a symmetric positive definite 
transmissibility matrix $T_K = (T_K^{\nu,\nu'})_{(\nu,\nu')\in \Xi_K\times \Xi_K}$ depending 
only on the cell $K$ geometry (including the choices of $\x_K$ and of 
$\x_\sigma, \sigma\in \faces_K$) and on the cell matrix diffusion tensor. 
The fracture fluxes are given by 
$$
F_{\sigma,\s}(u_\D) = \sum_{s\in \nodes_\sigma} T_\sigma^{\s,\s'} (u_\sigma - u_{s'}), 
$$
with a symmetric positive definite 
transmissibility matrix $T_\sigma = (T_\sigma^{\s,\s'})_{(\s,\s')\in \nodes_\sigma\times \nodes_\sigma}$ depending 
only on the fracture face $\sigma$ geometry (including the choice of $\x_\sigma$) 
and on the fracture face width and tangential diffusion tensor. 
Let us refer to \cite{BGGM14} for a more detailed presentation and for the definition 
of $T_K$ and $T_\sigma$. 
 
The construction of the control volumes at each degree of freedom is based on partitions 
of the cells and of the fracture faces. These partitions are respectively denoted, for all $K\in\cells$, by
$$\overline{K} ~ =  ~\overline{\omega}_K ~ \bigcup ~ \left( \bigcup_{\s\in\nodes_K\setminus\nodes_D}\overline{\omega}_{K,\s} 
\right), $$
and, for all $\sigma\in\faces_\G$, by 
$$\overline{\sigma}~=~\overline{\Sigma}_\sigma~\bigcup ~ \left( \bigcup_{\s\in\nodes_\sigma\setminus\nodes_D}\overline{\Sigma}_{\sigma,\s} 
\right). $$
The practical implementation of the scheme does not require to build explicitly the geometry of these partitions but only need to define the matrix volume fractions 
$$
\alpha_{K,\s} = {\int_{\omega_{K,\s}} d\x \over \int_K d\x}, 
\s\in\nodes_K\setminus(\nodes_D\cup \nodes_\Gamma), K\in\cells, 
$$
constrained to satisfy $\alpha_{K,\nu}\geq 0$,  
and $\sum_{\s\in\nodes_K\setminus(\nodes_D\cup\nodes_\Gamma)}\alpha_{K,\s} \leq 1$, 
as well as the fracture volume fractions 
$$
\alpha_{\sigma,\s} = {\int_{\Sigma_{\sigma,\s}} d_f(\x) d\tau(\x) \over \int_\sigma d_f(\x) d\tau(\x)}, 
\s\in\nodes_\sigma\setminus\nodes_D, \sigma\in\faces_\G, 
$$
constrained to satisfy 
$\alpha_{\sigma,\s}\geq 0$, and $\sum_{\s\in\nodes_\sigma\setminus\nodes_D}\alpha_{\sigma,\s}\leq 1$, where we denote by $d\tau({\bf x})$ the $2$ dimensional Lebesgue measure on $\Gamma$. 
Let us also set 
$$
\phi_K=   (1-\sum_{\s\in \nodes_K\setminus(\nodes_D\cup \nodes_\Gamma)}\alpha_{K,\s})\int_K \phi_m(\x) d\x
 \quad \mbox{ for } K\in \cells, 
$$ 
and 
$$
\phi_\sigma = (1-\sum_{\s\in \nodes_\sigma\setminus\nodes_D}\alpha_{\sigma,\s})
\int_\sigma \phi_f(\x) d_f(\x) d\tau(\x) 
\quad \mbox{ for } \sigma\in \faces_\Gamma, 
$$ 
as well as 
$$
\phi_{\s} =   \sum_{K\in \cells_\s} \alpha_{K,\s}\int_K \phi_m(\x) d\x \quad \mbox{ for } \s \in \nodes\setminus (\nodes_D\cup \nodes_\Gamma), 
$$
and
$$ 
\phi_{\s} = \sum_{\sigma \in \faces_{\Gamma,\s}}\alpha_{\sigma,\s}\int_\sigma \phi_f(\x) d_f(\x) d\tau(\x) 
\quad \mbox{ for } \s \in \nodes_\Gamma\setminus \nodes_D, 
$$
which correspond to the porous volumes distributed to the degrees of freedom excluding the  
Dirichlet nodes. 
The rock complementary volume in each control volume $\nu\in \cells\cup \faces_\Gamma\cup(\nodes\setminus\nodes_D)$ is denoted by $\bar \phi_\nu$.

As shown in \cite{BGGM14}, the flexibility in the choice of the control volumes is a crucial asset, compared with usual CVFE approaches and allows to significantly improve the accuracy of the scheme when the permeability field is highly heterogeneous. As exhibited in Figure \ref{fig_vag_CV}, as opposed to usual 
CVFE approaches, this flexibility allows to define the control volumes 
in the fractures with no contribution from the matrix in order to avoid to artificially enlarge the flow path in the fractures.  

A rocktype is assigned to each cell, node and fracture face.
In our case, for cells and for nodes not located along the fractures, the matrix rocktype is assigned. For fracture nodes and faces at the interface between the matrix and the fracture rocktypes, the fracture rocktype is assigned corresponding to the most pervious rock type consistently with the choice of the control volumes (see \cite{BGGM14}).
For convenience's sake, in the following, we will denote by $k_{r,\nu}^\alpha$ the corresponding relative permeability function for $\nu\in\cells\cup\nodes\cup\faces_\Gamma$. \\

In the following, we will keep the notation $F_{K,\s}$, $F_{K,\sigma}$, $F_{\sigma,\s}$ for 
the VAG Darcy fluxes defined with the cellwise constant 
matrix permeability $\K_m$ and the facewise constant fracture width $d_f$ and tangential 
permeability $\K_f$. Since the rock properties are fixed,
the VAG Darcy fluxes transmissibility matrices $T_K$ and $T_\sigma$ are computed only once. 

The VAG Fourier fluxes are denoted in the following 
by $G_{K,\s}$, $G_{K,\sigma}$, $G_{\sigma,\s}$. 
They are obtained with the isotropic matrix and fracture thermal conductivities 
averaged in each cell and 
in each fracture face using the previous time step fluid properties.
Hence VAG Fourier fluxes transmissibility matrices need 
to be recomputed at each time step. 

\subsection{Multi-branch non-isothermal well model}
\label{subsec_well}
Let $\mathcal{W}$ denote the set of wells. 
Each multi-branch well $\omega\in \mathcal{W}$ is defined by a set of oriented edges of the mesh
assumed to define a rooted tree oriented away from the root. This orientation corresponds to the drilling direction 
of the well.  The set of nodes of a well $\omega\in \mathcal{W}$ is denoted by $\nodes_\omega\subset\nodes$ and
its root node is denoted by $\s_\omega$.  A partial ordering is defined on the set of vertices $\mathcal{V}_\omega$
with $\s \underset{\omega}{<} \s'$ if and only if the unique path from the root $\s_\omega$ to $\s'$ passes through $\s$.
The set of edges of the well $\omega$ is denoted by $\mathcal{E}_\omega$ and for each edge $\welledge \in \edges_\omega$
we set $\welledge={\s}{\s'}$ with ${\s} \underset{\omega}{<} {\s'}$ (i.e. $\s$ is the parent node of $\s'$, see Figure \ref{fig_wellmodel}). It is assumed that
%$\edges_{\omega_1}\cap \edges_{\omega_2} = \emptyset$
$\nodes_{\omega_1}\cap \nodes_{\omega_2} = \emptyset$ for any $\omega_1, \omega_2 \in \mathcal{W}$ such that $\omega_1 \neq \omega_2$.

We focus on the part of the well that is connected to the reservoir through open hole, production liners or perforations. In this section, exchanges with the reservoir are dominated by convection and we decided to neglect heat losses as a first step. The latest shall be taken into account when modeling the wellbore flow up to the surface.
It is assumed that the radius $r_\omega$ of each well $\omega\in \mathcal{W}$ is small compared to the cell sizes in the neighborhood of the well. It results that the Darcy flux between the reservoir and the well at a
given well node $\s\in \mathcal{V}_\omega$ is obtained using the Two Point Flux Approximation
$$
V_\s^\omega = \WI_\s (p_\s - p_\s^\omega),  
$$
where $p_\s$ is the reservoir pressure at node $\s$ and $p_\s^\omega$ is the well pressure at node $\s$. The Well Index $\WI_{\s}$ is typically computed using Peaceman's approach (see \cite{Peaceman78,Peaceman83,CZ09}) and takes into account the unresolved singularity of the pressure solution in the neighborhood of the well. 
Fourier fluxes between the reservoir and the well could also be discretized using such Two Point Flux Approximation
but they are assumed to be small compared with thermal convective fluxes and will be neglected in the following well model. At each well node $\s\in \mathcal{V}_\omega$ the temperature inside the well is denoted by $T_\s^\omega$ and the volume fractions by $s_{\s,\omega}^\alpha$, $\alpha\in \{\l,\g\}$. The temperature in the reservoir at node $\s$ is denoted by $T_\s$, the saturations by $s_{\s}^\alpha$, and the phase mass fractions by $c^\alpha_\s$ for $\alpha\in \{\l,\g\}$.

For any $a\in \R$, let us define $a^+ = \max(a,0)$ and $a^- = \min(a,0)$. 
The mass flow rates between the reservoir and the well $\omega$ at a given node $\s \in \mathcal{V}_\omega$ are defined by the following phase based
upwind approximation of the mobilities:
\begin{equation}
  \label{eq_q1sw}
  \begin{array}{r@{\,\,}c@{\,\,}l}  
    q^{r\rightarrow \omega}_{\s,\alpha} &=&  \dsp \beta^{inj}_{\omega} {\rho^\alpha(p_\s^\omega,T_\s^\omega) \over \mu^\alpha(p_\s^\omega,T_\s^\omega)} k^\alpha_{r,\s}(s_{\s,\omega}^\alpha)(V_\s^\omega)^- +
    \beta^{prod}_{\omega}  c^\alpha_{\s}{\rho^\alpha(p_\s,T_\s) \over \mu^\alpha(p_\s,T_\s)} k^\alpha_{r,\s}(s_{\s}^\alpha) (V_\s^\omega)^+,\\[2ex]
    q^{r\rightarrow \omega}_{\s,\mass} &=& \dsp \sum_{\alpha\in \{\l,\g\}} q^{r\rightarrow \omega}_{\s,\alpha}, 
    \end{array}
\end{equation}
and the energy flow rate is defined similarly by 
\begin{equation}
\label{eq_q2sw}
q^{r\rightarrow \omega}_{\s,e} = \sum_{\alpha\in \{\l,\g\}} h^\alpha(p_\s^\omega,T_\s^\omega) (q^{r\rightarrow \omega}_{\s,\alpha})^- +
h^\alpha(p_\s,T_\s) (q^{r\rightarrow \omega}_{\s,\alpha})^+.
\end{equation}

The well coefficients $\beta^{inj}_{\omega}$ and  $\beta^{prod}_{\omega}$ are used to impose specific well behavior. The general case corresponds to $\beta^{inj}_{\omega} = \beta^{prod}_{\omega} = 1$.
Yet, for an injection well, it will be convenient as explained in subsection \ref{sec_inj_well}, to impose that the mass flow rates  $q^{r\rightarrow \omega}_{\s,\mass}$ are non positive for all nodes $s\in \mathcal{V}_\omega$
corresponding to set $\beta^{inj}_{\omega} = 1$ and $\beta^{prod}_{\omega} = 0$.
Likewise, for a production well, it will be convenient 
as explained in subsection \ref{sec_prod_well}, to set $\beta^{inj}_{\omega} = 0$ and $\beta^{prod}_{\omega} = 1$ which corresponds to assume that the mass flow rates $q^{r\rightarrow \omega}_{\s,\mass}$ are non negative for all nodes $s\in \mathcal{V}_\omega$. These simplifying options currently prevent the modeling of cross flows where injection and production occur in different places of the same well, as it sometimes happen in geothermal wells, typically in closed wells.

\subsubsection{Well physical model}
Our conceptual model inside the well assumes that the flow is stationary at the reservoir time scale
along with perfect mixing and thermal equilibrium. The Fourier fluxes and the wall friction are neglected and the pressure distribution is assumed hydrostatic along the well. 

For the sake of simplicity, 
the flow rate between the reservoir and the well is considered concentrated at each node $\s$ of the well.
For each edge $\welledge\in \mathcal{E}_\omega$, let us denote by $q^\alpha_\welledge$ the mass flow rate of phase $\alpha$ along the edge $\welledge$ oriented positively from ${\s'}$ to ${\s}$ with $\welledge={\s}{\s'}$ (let us recall that ${\s}$ is the parent node of ${\s'}$). 

\begin{figure}[htbp!]
\begin{center}
\includegraphics[width=0.35\textwidth]{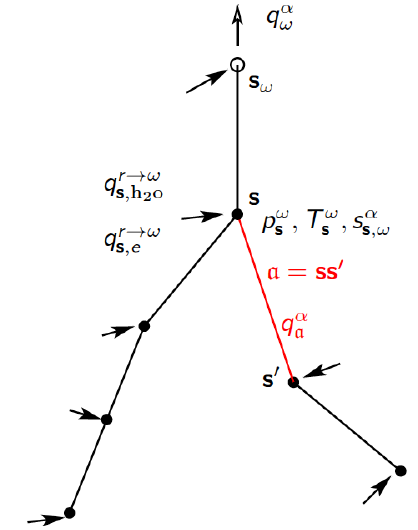}   
\caption{Example of multi-branch well $\omega$ with its root node $\s_\omega$, one edge $\mathfrak{a}=\s\s'$ and the main physical quantities: the well mass flow rates $q^\alpha_\omega$, the mass and energy flow rates between the reservoir and the well $q^{r\rightarrow \omega}_{\s,\mass}$, $q^{r\rightarrow \omega}_{\s,\energy}$, the well node pressure, temperature and saturations $p^\omega_\s, T^\omega_\s, s^\alpha_{\s,\omega}$, and the edge mass flow rates $q_{\mathfrak{a}}^\alpha$. }
\label{fig_wellmodel}
\end{center}
\end{figure}

Let $\alpha\in \{\l,\g\}$, the set of well unknowns is defined at each node
$\s\in \nodes_\omega$ by the well pressure $p_\s^\omega$, the well temperature $T_\s^\omega$, the well saturations $s_{\s,\omega}^\alpha$, and at each  edge $\welledge \in \edges_\omega$ by the  mass flow rates $q^\alpha_\welledge$. 
These well unknowns are complemented by the well mass flow rates $q^\alpha_\omega$ which are non negative for production wells and non positive for injection wells (see Figure \ref{fig_wellmodel}).  \\

For each edge $\welledge={\s}{\s'}\in \edges_\omega$, and each phase $\alpha$, let us define the following phase based upwind approximations of the specific enthalpy, mass density and saturation
\begin{equation}
  \label{eq_ha_rhoa_sa}
h^\alpha_\welledge =
\left\{\begin{array}{r@{\,\,}c@{\,\,}l}
& h^\alpha(p_{{\s'}}^\omega,T_{{\s'}}^\omega) \mbox{ if } q^\alpha_\welledge \geq 0,\\[1ex]
& h^\alpha(p_{{\s}}^\omega,T_{{\s}}^\omega) \mbox{ if } q^\alpha_\welledge < 0.  
\end{array}\right.
\rho^\alpha_\welledge =
\left\{\begin{array}{r@{\,\,}c@{\,\,}l}
& \rho^\alpha(p_{{\s'}}^\omega,T_{{\s'}}^\omega) \mbox{ if } q^\alpha_\welledge \geq 0,\\[1ex]
& \rho^\alpha(p_{{\s}}^\omega,T_{{\s}}^\omega) \mbox{ if } q^\alpha_\welledge < 0.  
\end{array}\right.
s^\alpha_\welledge =
\left\{\begin{array}{r@{\,\,}c@{\,\,}l}
& s^\alpha_{\s',\omega} \mbox{ if } q^\alpha_\welledge \geq 0,\\[1ex]
& s^\alpha_{\s,\omega} \mbox{ if } q^\alpha_\welledge < 0.  
\end{array}\right.
\end{equation}
For all ${\s}{\s'} = \welledge\in \edges_\omega$, let us set $\kappa_{\welledge,\s'} = -1$ and $\kappa_{\welledge,\s} = 1$.
The well equations account for the 
mass and energy conservations at each node of the well combined with the sum to one of the saturations and the thermodynamical equilibrium. Let $\edges_\s^\omega\subset \edges_\omega$ denote the set of well edges sharing the node $\s\in \nodes_\omega$, then 
for all $\s\in \mathcal{V}_\omega$ we obtain the equations 
\begin{equation}
  \label{well_cons1}
\left\{\begin{array}{r@{\,\,}c@{\,\,}l}
& \dsp q^{r\rightarrow \omega}_{\s,\mass} + \sum_{\welledge\in \edges_\s^\omega}  \sum_{\alpha\in \{\l,\g\}} \kappa_{\welledge,\s} q^\alpha_\welledge  =\delta_\s^{\s_\omega} \sum_{\alpha\in \{\l,\g\}} q^\alpha_\omega,\\[2ex]
& \dsp q^{r\rightarrow \omega}_{\s,e} + \sum_{\welledge\in \edges_\s^\omega}  \sum_{\alpha\in \{\l,\g\}} h^\alpha_\welledge \kappa_{\welledge,\s} q^\alpha_\welledge    = 
\delta_\s^{\s_\omega} \sum_{\alpha\in \{\l,\g\}} \left( \bar h^\alpha_\omega (q^\alpha_\omega)^- + h^\alpha(p_\s^\omega, T_\s^\omega) (q^\alpha_\omega)^+ \right), \\[2ex]
& s_{\s,\omega}^\l + s_{\s,\omega}^\g = 1,\\[2ex]
& p_\s^\omega = p_{\rm sat}(T_\s^\omega) \mbox{ if } s_{\s,\omega}^\g > 0 \mbox{ and } s_{\s,\omega}^\l > 0, \\[2ex]
& p_\s^\omega \geq p_{\rm sat}(T_\s^\omega) \mbox{ if } s_{\s,\omega}^\g = 0, \quad p_\s^\omega \leq p_{\rm sat}(T_\s^\omega) \mbox{ if } s_{\s,\omega}^\l=1,
\end{array}\right.
\end{equation}
where $\delta$ stands for the Kronecker symbol, and $\bar h^\alpha_\omega$ for prescribed specific enthalpies in the case of injection wells.
Inside the well, the hypothesis of hydrostatic pressure distribution implies that
\begin{equation}
  \label{well_hydrostat}
p_{{\s}}^\omega - p_{{\s'}}^\omega + \rho_\welledge g (z_{{\s}} - z_{{\s'}}) = 0, 
\end{equation}
for each edge ${\s}{\s'} = \welledge\in \edges_\omega$, where $\rho_\welledge$ is the mass density of the liquid gas mixture. The system is completed  by a slip closure law expressing the slip between the liquid velocity $u^\l_\welledge$  and the gas velocity $u^\g_\welledge$  at each edge $\welledge\in \edges_\omega$ with
$$
q^\alpha_\welledge = \pi r_\omega^2 \rho^\alpha_\welledge  s^\alpha_\welledge  u^\alpha_\welledge. 
$$
In the following simplified well models developed in subsections \ref{sec_inj_well} and \ref{sec_prod_well}, a zero slip law will be assumed for simplicity in such a way that $u^\l_\welledge = u^\g_\welledge$. Note that these simplified well models could be easily extended to account for non-zero slip laws as well as for an explicit approximation of the wall friction along the wells. The two fundamental assumptions to obtain these simplified well models are
\begin{itemize}
\item[(i)] prescribed sign of the mass flow rates $q^{r\rightarrow \omega}_{\s,\alpha}$, $\s\in \nodes_\omega$, forced to be all non-negative for production wells and all non-positive for injection wells, 
\item[(ii)] neglected Fourier fluxes compared with thermal convection fluxes. 
\end{itemize}

The well boundary conditions prescribe a limit total mass flow rate $\bar q_\omega$ and
a limit bottom hole pressure $\bar p_\omega$. Then, complementary constraints accounting for usual well monitoring conditions, are imposed 
between $q_\omega-\bar q_\omega$ and $p_\omega - \bar p_\omega$ using the notations 
$$
p_\omega = p_{\s_\omega}^\omega \mbox{ and } q_\omega = \sum_{\alpha\in \{\l,\g\}} q^\alpha_\omega. 
$$
In the following subsections, we consider the particular case of injection wells assuming a pure liquid phase, and the case of production wells. 
The flow rates are enforced to be non positive (resp. non negative)
at all well nodes for injection wells (resp. production wells).
It corresponds to set $\beta^{inj}_{\omega} = 1$,
$\beta^{prod}_{\omega} = 0$ for an injection well and
$\beta^{inj}_{\omega} = 0$, $\beta^{prod}_{\omega} = 1$ for a production well.
The limit bottom hole pressure $\bar p_\omega$ is a maximum (resp. minimum) pressure and the limit total mass flow rate $\bar q_\omega$ is
a minimum non positive (resp. maximum non negative) flow rate for injection (resp. production) wells.

In both cases, using an explicit computation of the hydrostatic pressure drop, the well model will be reduced
to a single equation and a single implicit unknown corresponding to the well reference pressure $p_\omega$
(see \textit{e.g.} \cite{Aunzo91}). 

\subsubsection{Liquid injection wells}
\label{sec_inj_well}
The injection well model sets $\beta^{inj}_{\omega} = 1$, $\beta^{prod}_{\omega} = 0$ and prescribes the minimum well total mass flow rate
$\bar q_\omega \leq 0$, the well maximum bottom hole pressure
$\bar p_\omega$ and the well specific liquid enthalpy $\bar h^\l_\omega$. It is assumed that the injection is in liquid phase and that no gas will appear in the well during the simulation as it is usually the case in geothermal systems. 

Since $\beta^{inj}_{\omega} = 1$ and $\beta^{prod}_{\omega} = 0$, the mass flow rates $q^\alpha_\welledge$ are enforced to be non negative and it results from
\eqref{well_cons1}, and the assumption that the gas phase does not appear in the well 
that $h^\l_\welledge = \bar h^\l_\omega$ for all $\welledge\in \edges_\omega$ and that $s_{\s,\omega}^\l =  1-  s_{\s,\omega}^\g = 1$ for all $\s\in \nodes_\omega$.

Given the previous time step well reference pressure $p_\omega^{n-1} = p_{\s_\omega}^{\omega,n-1}$, we first compute the pressures along the well solving the equations
\begin{align*}
  & p_{\s}^\omega - p_{\s'}^\omega + \rho_\welledge g (z_{\s} - z_{\s'}) = 0 \mbox{ for all } \welledge=\s\s'\in \edges_\omega,\\[1ex]
  & p_{\s_\omega}^\omega = p_{\s_\omega}^{\omega,n-1},\\[1ex]
  & \rho_\welledge = \rho^\l(p_{\s}^\omega,T_{\s}^\omega) \mbox{ for all } \welledge=\s\s'\in \edges_\omega,\\[1ex]
  & h^\l(p^{\omega}_{\s},T^{\omega}_{\s}) = \bar h^\l_\omega \mbox{ for all } \s \in \nodes_\omega. 
\end{align*}
We deduce the explicit pressure drops
$$
\Delta p^{\omega,n-1}_\s = p_{\s}^\omega - p_\omega^{n-1},  
$$
which provide for all $\s\in \nodes_\omega$ the pressures $p^{\omega,n}_\s$ and temperatures $T^{\omega,n}_{\s}$ along the well at the current time step $n$ such that 
\begin{align*}
  & p^{\omega,n}_\s =  p^n_{\omega} + \Delta p^{\omega,n-1}_\s,\\[1ex]
  & h^\l(p^{\omega,n}_{\s},T^{\omega,n}_{\s}) = \bar h^\l_\omega.  
\end{align*}

The mass and energy flow rates at each node $\s\in \mathcal{V}_\omega$ between the reservoir and the well are defined
by \eqref{eq_q1sw}-\eqref{eq_q2sw} with $\beta^{inj}_{\omega} = 1$ and $\beta^{prod}_{\omega} = 0$ and depend
only on the implicit unknowns $p^n_\omega$ and $p_\s^n$. They are respectively denoted by
$q_{\s,\mass}^{r\rightarrow \omega}(p^n_\s, p^{n}_\omega)$ and $q_{\s,e}^{r\rightarrow \omega}(p^n_\s, p^{n}_\omega)$. 

The well equation at the current time step is defined by the following complementary constraints  between
the prescribed minimum well total mass flow rate and the prescribed maximum bottom hole pressure
%$$
%\min(\sum_{\s\in \nodes_{\omega} } q_{\s,\mass}^{r\rightarrow \omega}(p^n_\s,T_\s^n, p^{\omega,n}_\s)    - \bar q_\omega, \bar p_\omega-p^n_\omega) = 0. 
%$$
\begin{equation}
  \label{injection_well_constraints}
  \left\{\begin{array}{r@{\,\,}c@{\,\,}l}
 \dsp \(\sum_{\s\in \nodes_{\omega} }  q_{\s,\mass}^{r\rightarrow \omega}(p^n_\s, p^{n}_\omega)  - \bar q_\omega\) \(\bar p_\omega-p^n_\omega\) &=& 0, \\
 \dsp \sum_{\s\in \nodes_{\omega} } q_{\s,\mass}^{r\rightarrow \omega}(p^n_\s, p^{n}_\omega)   - \bar q_\omega &\geq& 0, \\
 \dsp \bar p_\omega-p^n_\omega &\geq& 0. 
\end{array}\right.
  \end{equation}

\subsubsection{Production wells}
\label{sec_prod_well}
The production well model sets $\beta^{inj}_{\omega} = 0$, $\beta^{prod}_{\omega} = 1$ and prescribes the maximum well total mass flow rate $\bar q_\omega \geq 0$ and the well minimum bottom hole pressure $\bar p_\omega$.

The solution at the previous time step $n-1$ provides the pressure drop $\Delta p^{\omega,n-1}_\s$ at each node $\s\in \nodes_\omega$.  This computation based on thermodynamical equilibrium is detailed below. 
As for the injection well, we deduce the well pressures using the bottom well pressure at the current time step $n$ 
$$
p^{\omega,n}_\s =  p^n_{\omega} + \Delta p^{\omega,n-1}_\s. 
$$
The mass and energy flow rates at each node $\s\in \mathcal{V}_\omega$ between the reservoir and the well
are defined by \eqref{eq_q1sw}-\eqref{eq_q2sw} with $\beta^{inj}_{\omega} = 0$ and $\beta^{prod}_{\omega} = 1$ and depend
only on the implicit reservoir unknowns $X^n_\s$ setting
$$
X_\s = \(P_\s,T_\s, s^\l_\s, s^\g_\s, c^\l_\s, c^\g_\s\),   
$$
and on the implicit well unknown $p^n_\omega$.
They are respectively denoted by $q_{\s,\mass}^{r\rightarrow \omega}(X^n_\s,p^{n}_\omega)$ and $q_{\s,e}^{r\rightarrow \omega}(X^n_\s, p^{n}_\omega)$.

The well equation at the current time step is defined by the following complementary constraints  between
the prescribed maximum well total mass flow rate and the prescribed minimum bottom hole pressure
\begin{equation}
  \label{production_well_constraints}
   \left\{\begin{array}{r@{\,\,}c@{\,\,}l}
\dsp \( \bar q_\omega - \sum_{\s\in \nodes_{\omega} } q_{\s,\mass}^{r\rightarrow \omega}(X^n_\s,p^{n}_\omega) \) \(p_\omega^{n}-\bar p_\omega\) &=& 0, \\
\dsp \bar q_\omega - \sum_{\s\in \nodes_{\omega} } q_{\s,\mass}^{r\rightarrow \omega}(X^n_\s,p^{n}_\omega) &\geq & 0,\\
\dsp p_\omega^{n}-\bar p_\omega &\geq & 0.
\end{array}\right.
\end{equation}

Let us now detail the computation of the pressure drop at each node $\s\in \mathcal{V}_\omega$ using the previous time step
solution $n-1$ consisting of the reservoir unknowns and the well pressures. 
We first compute the well temperature $T_{\s}^{\omega,n-1}$ and saturations $s^{\alpha,n-1}_{\s,\omega}$ at each node $\s$ using equations \eqref{well_cons1}. Summing the mass and energy equations of \eqref{well_cons1}  over all nodes $\s'' \underset{\omega}{\geq} \s$, we obtain for all $\welledge=\s'\s\in \edges_\omega$ that

\begin{align*}
    &   \sum_{\alpha\in \{\l,\g\}}  Q^{\alpha,n-1}_\welledge =  \sum_{\s''\in {\cal V}_\omega | \s'' \underset{\omega}{\geq} \s} q_{\s'',\mass}^{r\rightarrow \omega}(X^{n-1}_{\s''},p^{n-1}_{\omega}) = Q_{\s,\mass}^\omega,\\
      &   \sum_{\alpha\in \{\l,\g\}}  h^\alpha(p_\s^{\omega,n-1},T_\s^{\omega,n-1}) Q^{\alpha,n-1}_\welledge  =   \sum_{\s'\in {\cal V}_\omega | \s'' \underset{\omega}{\geq} \s} q_{\s',e}^{r\rightarrow \omega}(X^{n-1}_{\s''},p^{n-1}_{\omega}) = Q_{\s,e}^\omega,  
\end{align*}
with
$$
Q^{\alpha,n-1}_\welledge  = \pi r_\omega^2  \rho^\alpha(p_\s^{\omega,n-1},T_\s^{\omega,n-1}) s_{\s,\omega}^{\alpha,n-1}  u^{\alpha,n-1}_\welledge, \, \alpha\in \{\l,\g\}. 
$$

It results that the thermodynamical equilibrium at fixed well pressure $p^{\omega,n-1}_\s$, mass $Q_{\s,\mass}^\omega$ and energy $Q_{\s,e}^\omega$ provides the 
well temperature $T_\s^{\omega,n-1}$ and the well saturations $s_{\s,\omega}^{\alpha,n-1}$ at node $\s$ as follows.
Let us set $p=p^{\omega,n-1}_\s$. We first assume that both phases are present which implies that $T_{\rm sat} = (p_{\rm sat})^{-1}(p)$ and that the liquid mass fraction is given by   
$$
c^\l = { h^\g(p,T_{\rm sat}) - {Q_{\s,e}^\omega \over Q_{\s,\mass}^\omega}  \over h^\g(p,T_{\rm sat}) - h^\l(p,T_{\rm sat})}. 
$$
The following alternatives are checked: 
\begin{itemize}
\item[] {\it Two-phase state}: if $0 < c^\l <1$, the two-phase state is confirmed. Using the zero slip assumption, we obtain 
$$
\dsp T_\s^{\omega,n-1} = T_{\rm sat} \mbox{ and } 
s_{\s,\omega}^{\l,n-1} = 1-s_{\s,\omega}^{\g,n-1} = { { c^\l \over \rho^\l(p,T_{\rm sat})} \over { c^\l \over \rho^\l(p,T_{\rm sat})} + { 1-c^\l \over \rho^\g(p,T_{\rm sat})} }. 
$$
\item[] {\it Liquid state}: if $c^\l \geq 1$, then only the liquid phase is present, we set $s_{\s,\omega}^{\l,n-1} = 1$, $s_{\s,\omega}^{\g,n-1} = 0$, and $T_\s^{\omega,n-1}$ is the solution $T$ of
  $$
  h^\l(p,T) = {Q_{\s,e}^\omega \over Q_{\s,\mass}^\omega}.
  $$
\item[] {\it Gas state}: if $c^\l \leq 0$, then only the gas phase is present, we set $s_{\s,\omega}^{\l,n-1} = 0$, $s_{\s,\omega}^{\g,n-1} = 1$, and $T_\s^{\omega,n-1}$ is the solution $T$ of 
  $$
  h^\g(p,T) = {Q_{\s,e}^\omega \over Q_{\s,\mass}^\omega}.
  $$
\end{itemize}

Then, the explicit pressure drop
$$
\Delta p^{\omega,n-1}_\s = p_{\s}^\omega - p_\omega^{n-1},   
$$
is obtained from 
\begin{align*}
  & p_{\s}^\omega - p_{\s'}^\omega + \rho_\welledge g (z_{\s} - z_{\s'}) = 0 \mbox{ for all } \welledge=\s\s'\in \edges_\omega,\\[1ex]
  & p_{\s_\omega}^\omega = p_{\s_\omega}^{\omega,n-1},\\[1ex]
  & \rho_\welledge = \sum_{\alpha\in \{\l,\g\}} s_{\s,\omega}^{\alpha,n-1} \rho^\alpha(p_{\s}^{\omega,n-1},T_{\s}^{\omega,n-1}) \mbox{ for all } \welledge=\s\s'\in \edges_\omega. 
\end{align*}

\subsection{Discretization of the hybrid-dimensional non-isothermal two-phase flow model}

The time integration is based on a fully implicit Euler scheme to avoid severe restrictions on the time 
steps due to the small volumes and high velocities in the fractures. 
A phase based upwind scheme is used for the approximation of the 
mobilities in the mass and energy fluxes (see e.g. \cite{aziz-settari-79}). 
At the matrix fracture interfaces, we avoid mixing matrix 
and fracture rocktypes by choosing appropriate control volumes for $\sigma\in \faces_\Gamma$ 
and $\s\in \nodes_\Gamma$ (see Figure \ref{fig_vag_CV}). 
In order to avoid tiny control volumes at the  
nodes $\s\in \nodes_\Sigma$ located at the fracture intersection, 
the volume is distributed to such a node 
$\s$ from all the fracture faces containing the node $\s$. \\

For each $\nu \in \cells\cup \faces_\Gamma \cup \nodes$ the set of reservoir pressure, temperature, saturations and mass fractions unknowns 
is denoted by 
$
X_\nu = \(P_\nu,T_\nu, s^\l_\nu, s^\g_\nu, c^\l_\nu, c^\g_\nu\),   
$
where $c^\alpha_\nu$ is the mass fraction of the water component in phase $\alpha$ used to express the thermodynamical equilibrium. 
We denote by $X_\D$, the set of reservoir unknowns 
$$
X_\D = \{X_\nu, \, \nu \in \cells\cup \faces_\Gamma \cup \nodes\},  
$$
and similarly by $P_\D$ and $T_\D$ the sets of reservoir pressures and temperatures.   
The set of well bottom hole pressures is denoted by $P_\mathcal{W} = \{ p_\omega, \, \omega \in \mathcal{W}\}$.

The Darcy fluxes taking into account the gravity term are defined by 
\begin{equation}\label{flux_phase_vag}
\left\{\begin{array}{r@{\,\,}c@{\,\,}ll}
&V^\alpha_{K,\nu}(X_\D) &=& \dsp F_{K,\nu}(P_\D) - {\rho^\alpha(p_K,T_K) + \rho^\alpha(p_{\nu},T_\nu) \over 2 } F_{K,\nu}({\cal G}_\D), \quad \nu\in \Xi_K, K\in \cells,\\[1ex]
&V^\alpha_{\sigma,\s}(X_\D) &=& \dsp F_{\sigma,\s}(P_\D) - {\rho^\alpha(p_\sigma,T_\sigma) + \rho^\alpha(p_{\s},T_\s) \over 2 } F_{\sigma,\s}({\cal G}_\D), \quad \s\in \nodes_\sigma, \sigma \in \faces_\Gamma, 
\end{array}\right.
\end{equation}
where ${\cal G}_\D$ denotes the vector $({\bf g}\cdot {\bf x}_\nu)_{\nu \in \cells\cup \faces_\Gamma \cup \nodes}$.

For each Darcy flux, let us define the upwind control volume $cv^\alpha_{\mu,\nu}$ such that  
\begin{equation*}
cv^\alpha_{K,\nu} = 
\left\{\begin{array}{r@{\,\,}c@{\,\,}l}
K &\text{ if } & V^\alpha_{K,{\nu}}(X_\D) \geqslant 0 \\[1ex]
\nu & \text{ if } & V^\alpha_{K,{\nu}}(X_\D) < 0
\end{array}\right. 
\text{ for } \nu\in \Xi_K, K\in \cells,
\end{equation*}
for the matrix fluxes, and such that  
\begin{equation*}
cv^\alpha_{\sigma,\s} = 
\left\{\begin{array}{r@{\,\,}c@{\,\,}l}
\sigma &\text{ if } &V^\alpha_{\sigma,\s}(X_\D) \geqslant 0 \\[1ex]
\s &\text{ if } &V^\alpha_{\sigma,\s}(X_\D) < 0
\end{array}\right.
\text{ for }  \s\in \nodes_\sigma, \sigma \in \faces_\Gamma,
\end{equation*}
for fracture fluxes. Using this upwinding, the mass and energy fluxes are given by 
\begin{align*}
& q^\alpha_{\nu,\nu'}(X_\D) = 
  c^\alpha_{cv^\alpha_{\nu,\nu'}} {\rho^\alpha(p_{cv^\alpha_{\nu,\nu'}},T_{cv^\alpha_{\nu,\nu'}}) \over \mu^\alpha(p_{cv^\alpha_{\nu,\nu'}},T_{cv^\alpha_{\nu,\nu'}})}  k_{r,cv^\alpha_{\nu,\nu'}}^\alpha(s^\alpha_{cv^\alpha_{\nu,\nu'}})V^\alpha_{\nu,\nu'}(X_\D), \\
  & q^\mass_{\nu,\nu'}(X_\D) = \sum_{\alpha\in\{\l,\g\}} q^\alpha_{\nu,\nu'}(X_\D)
, \\
& q^\energy_{\nu,\nu'}(X_\D) = \sum_{\alpha\in\{\l,\g\}} h^\alpha(p_{cv^\alpha_{\nu,\nu'}},T_{cv^\alpha_{\nu,\nu'}} ) q^\alpha_{\nu,\nu'}(X_\D) + G_{\nu,\nu'}(T_\D). 
\end{align*}
In each control volume $\nu \in \cells\cup \faces_\Gamma \cup \nodes$, the mass and energy accumulations are denoted by   
\begin{align*}
  & {\cal A}_{\alpha,\nu}(X_\nu) =  \phi_\nu  \rho^\alpha(p_\nu,T_\nu) s^\alpha_\nu c^\alpha_\nu, \\
   & {\cal A}_{\mass,\nu}(X_\nu) = \sum_{\alpha\in\{\l,\g\}} {\cal A}_{\alpha,\nu}(X_\nu),\\
& {\cal A}_{\energy,\nu}(X_\nu) = \sum_{\alpha\in\{\l,\g\}} e^\alpha(p_\nu,T_\nu) {\cal A}_{\alpha,\nu}(X_\nu)  + \bar \phi_\nu E_r(p_\nu,T_\nu). 
\end{align*}
We can now state the system of discrete equations at each 
time step $n=1,\cdots,N_{t_f}$ which accounts for the mass ($i=\mass$) and energy ($i=\energy$) 
conservation equations in each cell $K\in \cells$: 
\begin{equation}
\label{VAGconvcell}
R_{K,i} (X_\D^{n}) 
:= \frac{ {\cal A}_{i,K}(X_K^{n})-{\cal A}_{i,K}(X_K^{n-1})}{\Delta t^n} 
+ \sum_{\s \in \nodes_K}  q^i_{K,\s}(X_\D^n)
+ \sum_{\sigma \in \faces_{\Gamma}\cap\faces_K}  q^i_{K,\sigma}(X_\D^n) = 0,
\end{equation}
in each fracture face $\sigma \in \faces_{\Gamma}$: 
\begin{equation}
\label{VAGconvfrac} 
R_{\sigma,i} (X_\D^{n}) :=
\frac{{\cal A}_{i,\sigma}(X_\sigma^{n})-{\cal A}_{i,\sigma}(X_\sigma^{n-1})}{\Delta t^n} + 
\sum_{\s \in {\cal V}_\sigma}  q^i_{\sigma,\s}(X_\D^n) 
+ \sum_{K \in {\cal M}_\sigma}  - q^i_{K,\sigma}(X_\D^n) = 0, 
\end{equation}
and at each node $\s\in \nodes\setminus\nodes_D$:
\begin{equation}
  \label{VAGconvnode}
  \begin{array}{r@{\,\,}c@{\,\,}l}
R_{\s,i} (X_\D^{n}, P_\mathcal{W}^n) &:=&\dsp 
\frac{{\cal A}_{i,\s}(X_\s^{n})-{\cal A}_{i,\s}(X_\s^{n-1})}{\Delta t^n} + 
\sum_{\sigma \in {\cal F}_{\Gamma,\s}}  - q^i_{\sigma,\s} (X_\D^n)
+ \sum_{K \in {\cal M}_\s}  - q^i_{K,\s}(X_\D^n)  \\
&& \qquad\qquad\qquad\qquad\qquad \dsp + \sum_{\omega\in \mathcal{W} | \s \in \mathcal{V}_\omega} q_{\s,i}^{r\rightarrow \omega}(X^n_\s,p^{\omega,n}_\s) = 0.
\end{array}
\end{equation}
It is coupled with the well equations for the injection wells $\omega\in \mathcal{W}_{inj}$ 
\begin{equation}
\label{injwelleq} 
R_{\omega} (X_\D^{n}, P_\mathcal{W}^n) := - \min(\sum_{\s\in {\cal V}_{\omega} } q_{\s,\mass}^{r\rightarrow \omega}(X^n_\s,p^{n}_\omega)   - \bar q_\omega, \bar p_\omega-p^n_\omega) = 0, 
\end{equation}
and for the production wells $\omega\in \mathcal{W}_{prod}$ 
\begin{equation}
\label{prodwelleq} 
R_{\omega} (X_\D^{n}, P_\mathcal{W}^n) :=  \min( \bar q_\omega - \sum_{\s\in {\cal V}_{\omega} } q_{\s,\mass}^{r\rightarrow \omega}(X^n_\s,p^{n}_\omega), p_\omega^{n}-\bar p_\omega) = 0, 
\end{equation}
reformulating respectively \eqref{injection_well_constraints} and \eqref{production_well_constraints} using the min function.

The system is closed with thermodynamical equilibrium and the sum to one of the saturations
\begin{equation}
  \label{closure_eq}
\begin{array}{r@{\,\,}c@{\,\,}lll}  
R_1 (X_\nu^{n}) &:=&c_\nu^{\g,n} p^n_\nu - p_{\rm sat}(T^n_\nu) c^{\l,n}_\nu &=& 0,\\ [1ex] 
R_2 (X_\nu^{n}) &:=&  \min(s^{\l,n}_\nu, 1-c^{\l,n}_\nu ) &=& 0, \\[1ex]
R_3 (X_\nu^{n}) &:=&  \min(s^{\g,n}_\nu, 1-c^{\g,n}_\nu ) &=& 0, \\[1ex]
R_4(X_\nu^{n}) &:=& s^{\l,n}_\nu + s^{\g,n}_\nu-1 &=& 0,
\end{array}
\end{equation}
at all control volumes $\nu \in \cells\cup \faces_\Gamma \cup \nodes\setminus \nodes_\D$ as well as the Dirichlet boundary conditions 
$$
X_\s^n = X_{\s,D}, 
$$
for all $\s\in \nodes_\D$. \\

Let us denote by $R_\nu$ the vector $\(R_{\nu,i}, \ i \in \{\mass, \energy\},\, R_j(X_\nu), j\in \{1,\cdots,4\}\) $, and 
let us rewrite the conservation and closure equations \eqref{VAGconvcell}, \eqref{VAGconvfrac}, 
\eqref{VAGconvnode}, \eqref{injwelleq}, \eqref{prodwelleq}, \eqref{closure_eq}   as well 
as the Dirichlet boundary conditions in vector form defining the following 
nonlinear system at each time step $n=1,2,...,N_{t_f}$
\begin{eqnarray}
\label{NonLinearSystem}
\mathbf{0} = \mathcal{R} (X_\D, P_\mathcal{W}) := 
\left\{\begin{array}{llllll}
R_{\s} (X_\D, P_\mathcal{W}), 
\,\, \s \in \nodes, \\
R_{\sigma} (X_\D) ,
\,\, \sigma \in \faces_\Gamma, \\
R_{K} (X_\D),
\,\, K \in \cells, \\
R_{\omega} (X_\D, P_\mathcal{W}), \,\, \omega \in \mathcal{W}, \\
\end{array}\right.
\end{eqnarray}
where the superscript $n$ is dropped to simplify the notations 
and where the Dirichlet boundary conditions 
have been included at each Dirichlet node $\s\in {\cal V}_D$ in order to 
obtain a system size independent of the boundary conditions. 

The nonlinear system ${\cal R}(X_\D, P_\mathcal{W})=0$ 
is solved by a Newton-min algorithm \cite{Kr11}. Our implementation is based on an active set method both for the well equations and the thermodynamical equilibrium.

For the well equations, we enforce either the total mass flow rate or the bottom hole pressure at each Newton iterate and use the remaining inequality constraint
to switch from prescribed total mass flow rate to prescribed bottom hole pressure and vice versa. 

For the thermodynamical equilibrium, we distinguish a two-phase state $I^n_\nu = \{\l,\g\}$, a liquid state $I^n_\nu = \{\l\}$, and a gas state $I^n_\nu = \{\g\}$.
For $I^n_\nu = \{\l,\g\}$, the closure equations provide $c^{\l,n}_\nu = c^{\g,n}_\nu = 1$,
$p^n_\nu - p_{\rm sat}(T^n_\nu)$ and $s^{\l,n}_\nu = 1- s^{\g,n}_\nu$ and we define $Y_\nu = (p^n_\nu,s^{\g,n}_\nu)$ as primary unknowns.
For $I^n_\nu = \{\l\}$, the closure equations provide $c^{\l,n}_\nu =  1$,
$c^{\g,n}_\nu = {p_{\rm sat}(T^n_\nu) \over p^n_\nu }$, $s^{\l,n}_\nu = 1$, $s^{\g,n}_\nu = 0$ and we define $Y_\nu = (p^n_\nu,T^{n}_\nu)$ as primary unknowns.
For $I^n_\nu = \{\g\}$, the closure equations provide $c^{\g,n}_\nu =  1$,
$c^{\l,n}_\nu = { p^n_\nu \over p_{\rm sat}(T^n_\nu)}$, $s^{\l,n}_\nu = 0$, $s^{\g,n}_\nu = 1$ and we define $Y_\nu = (p^n_\nu,T^{n}_\nu)$ as primary unknowns.
The inequality constraints are then used to switch from two-phase state to a one phase state and vice versa. \\

The Jacobian system at each Newton-min iteration is assembled w.r.t. the primary unknowns $Y_\D,P_\mathcal{W}$ and the mass and energy conservation equations \eqref{VAGconvcell}, \eqref{VAGconvfrac}, 
\eqref{VAGconvnode}, \eqref{injwelleq}, \eqref{prodwelleq}. The cell unknowns are locally eliminated without any additional fill-in before solving the linear system using 
the GMRES iterative solver preconditioned by a CPR-AMG preconditioner introduced in \cite{LVW:2001,SMW:2003}.
This preconditioner combines multiplicatively a parallel algebraic multigrid preconditioner (AMG) \cite{HY:2001} for a pressure block
of the linear system  with a block Jacobi ILU0 preconditioner for the full system.
In our case, the columns of the pressure block are defined by the node, the fracture face and the well pressure unknowns,
and its lines by the node and the fracture face mass conservation equations as well as the well equations. 

The parallel implementation is 
described in \cite{Xing.ea:2017} and \cite{Xing2018}. Let us recall that the distribution of wells to each MPI process $p$
is such that any well with a node belonging to the set of own nodes of $p$ belongs to the set of own and ghost wells of $p$. Then, the set of own and ghost nodes of $p$
is extended to include all the nodes belonging to
the own and ghost wells of $p$. These definitions ensure that (i) the local linearized systems can be assembled locally on each process without communication as in \cite{Xing.ea:2017}, 
and (ii) the pressure drops of the wells can be computed locally on each process without communication.
This last property is convenient since the pressure drop is a sequential computation along the well rooted tree. 
This parallelization strategy of the well model is based on the assumption that the number of additional ghost nodes resulting from the connectivity of the wells remains very small compared with the number of own nodes of the process.

\section{Numerical results}
\label{sec_num}

\subsection{Numerical convergence for a diphasic vertical well in an homogeneous reservoir}
\label{subsec_num_dipwell}

Let us consider the geothermal reservoir defined by the domain $\Omega = (-H, H)^2\times(0,H_z)$
where $H=1000$ m and $H_z=200$ m, including one vertical producer well 
along the line $\{(x, y, z) \in \Omega \,|\, x = y = 0\}$ of radius $r_{\omega}=0.1$ m.
The reservoir is assumed homogeneous with isotropic permeability 
$\K_m = k_m I, k_m = 5\times 10^{-14} \text{ m}^2$ and porosity $\phi_m = 0.15$.
It is assumed to be initially saturated with pure water in liquid phase. 
The  enthalpy, internal energy, mass density and viscosity of water in the liquid and gas phases are given from \cite{refthermo} by analytical laws as functions of the pressure and temperature.
The vapour pressure $P_{sat}(T)$ is given in Pa by the Clausius-Clapeyron equation 
$$
p_{\rm sat}(T) = 100 \exp\left(46.784 - \frac{6435}{T} - 3.868 \,\, log(T)\right). 
$$
The thermal conductivity is fixed to $\lambda_m = 2 \;\text{W}.\text{m}^{-1}.\text{K}^{-1}$,
 and the rock volumetric heat capacity is given by $C_r = 1.6\; \text{MJ}.\text{K}^{-1}.\text{m}^{-3}$ with $E_r(p,T) = C_r T$.  
The relative permeabilities are set to
$k_{r,m}^\alpha(s^\alpha) = (s^\alpha)^2$ for both phases $\alpha\in \{\l,\g\}$. 
The gravity vector is as usual ${\bf g}=(0,0,-g_z)$ with $g_z=9.81 \text{ m}.\text{s}^{-2}$.

The simulation consists in two stages both run on a family of refined uniform Cartesian meshes of size $n_x\times n_y \times n_z$ of the domain $\Omega$ with 
$$
(n_x,n_y,n_z) \in \{ (10,10,5),(20,20,10),(40,40,20),(80,80,40)\}. 
$$
These meshes  are labeled as  $\{h_1,\dots, h_4\}$ respectively. The well indexes are computed at each node of the well following \cite{Xing2018}.

At the first stage, the well is closed and a Dirichlet boundary condition
is imposed at the top
of the domain prescribing the pressure and the temperature equal to $p_m=4$ MPa
and $T_m = (p_{\rm sat})^{-1}(p_m)- 1\;\text{K}$; respectively,
and homogeneous Neumann boundary conditions are set at the bottom and at the sides of the domain.
This stage is run until the simulation reaches a stationary state with the liquid phase only, a constant temperature and an hydrostatic pressure depending only on the vertical coordinate. 

For the second stage, homogeneous Neumann boundary conditions are prescribed at the bottom and at the top of the domain $\Omega$,
but Dirichlet boundary conditions for the pressure and temperature are fixed at the sides of the domain to the ones at the end of stage one. 
The well is  set in an open state, i.e.,
it can produce,
and its monitoring conditions are defined by the minimum bottom hole pressure
$\bar{p}_\omega= 1$ bar (never reached in practice) and the maximum total mass flow rate $\bar{q}_\omega= 200  \text{ ton}.\text{hour}^{-1}$. The second stage is run on the time interval  $(0,t_f)$ with $t_f=30$ days.

\begin{figure}[H]
\begin{center}
\includegraphics[width=0.55\textwidth]{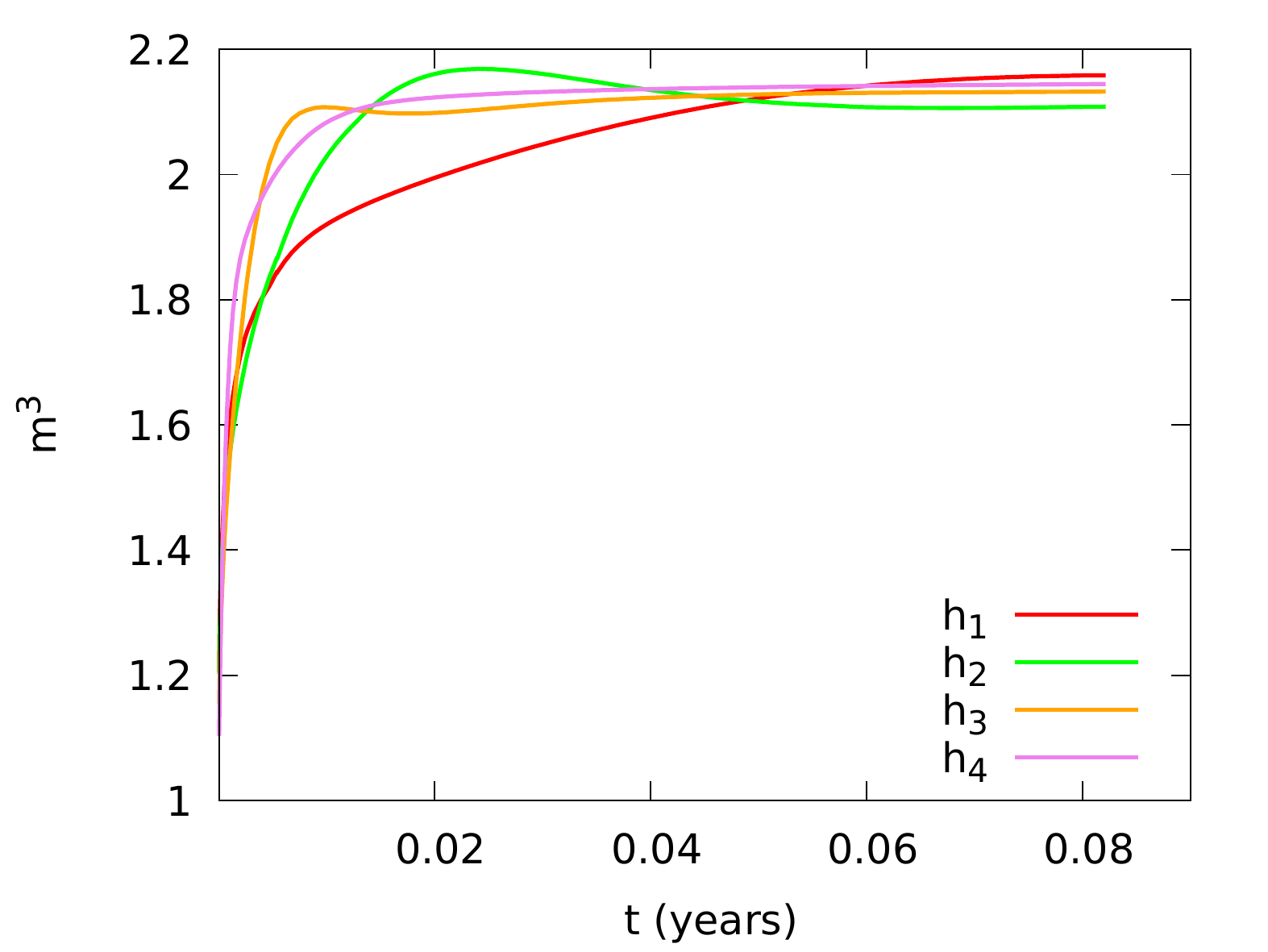}
    \caption{Total gas volume inside the well as a function of time on the different meshes.}
    \label{fig_gas_vol_well}
\end{center}
\end{figure}

\begin{figure}[H]
\begin{center}
\includegraphics[width=0.55\textwidth]{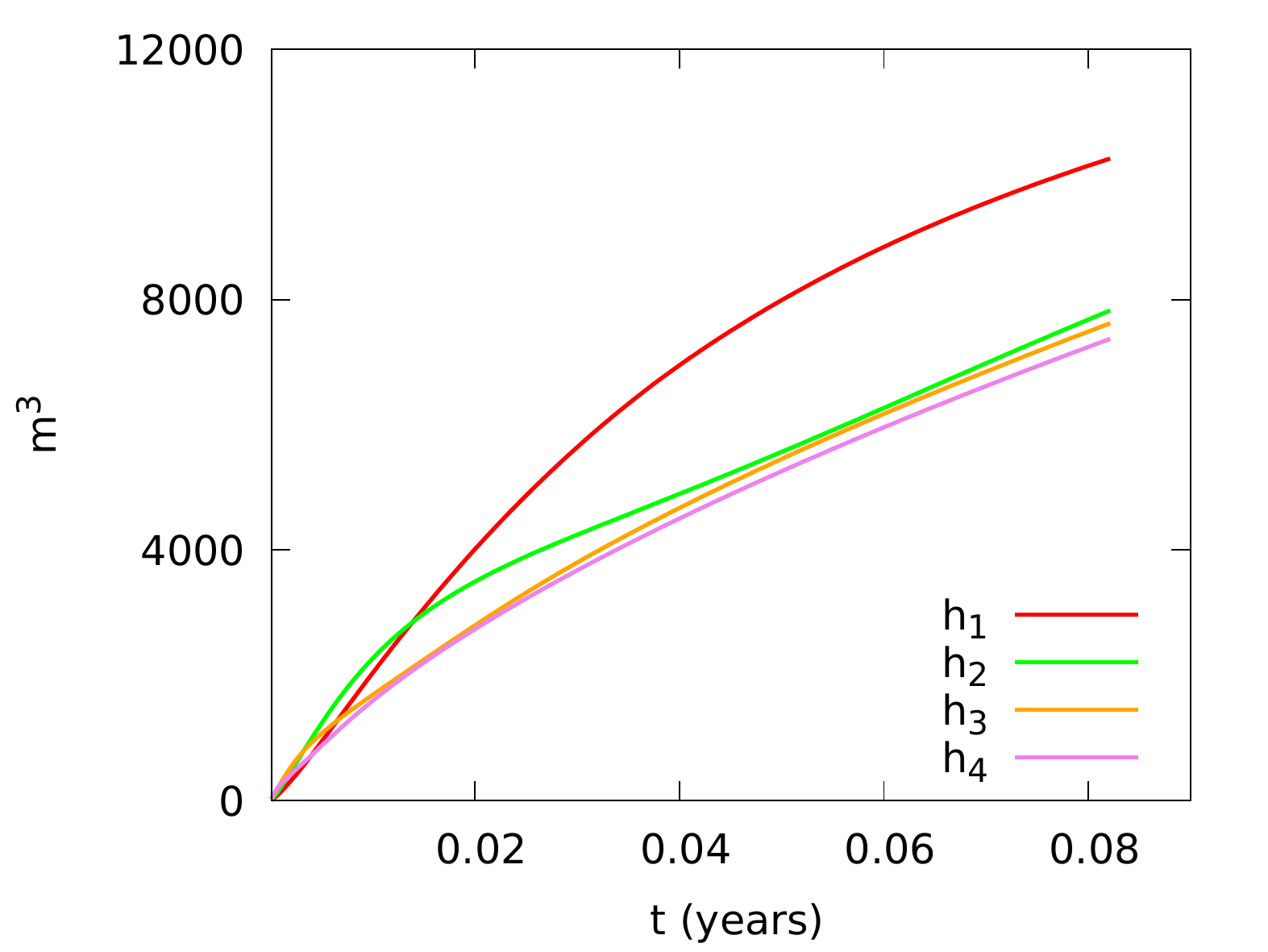}
    \caption{Total gas volume inside the reservoir  as a function of time on the different meshes.}
\label{fig_gas_vol_res}
\end{center}
\end{figure}

\begin{figure}[H]
    \centering
    \begin{subfigure}[b]{0.45\textwidth}
      \includegraphics[width=\textwidth]{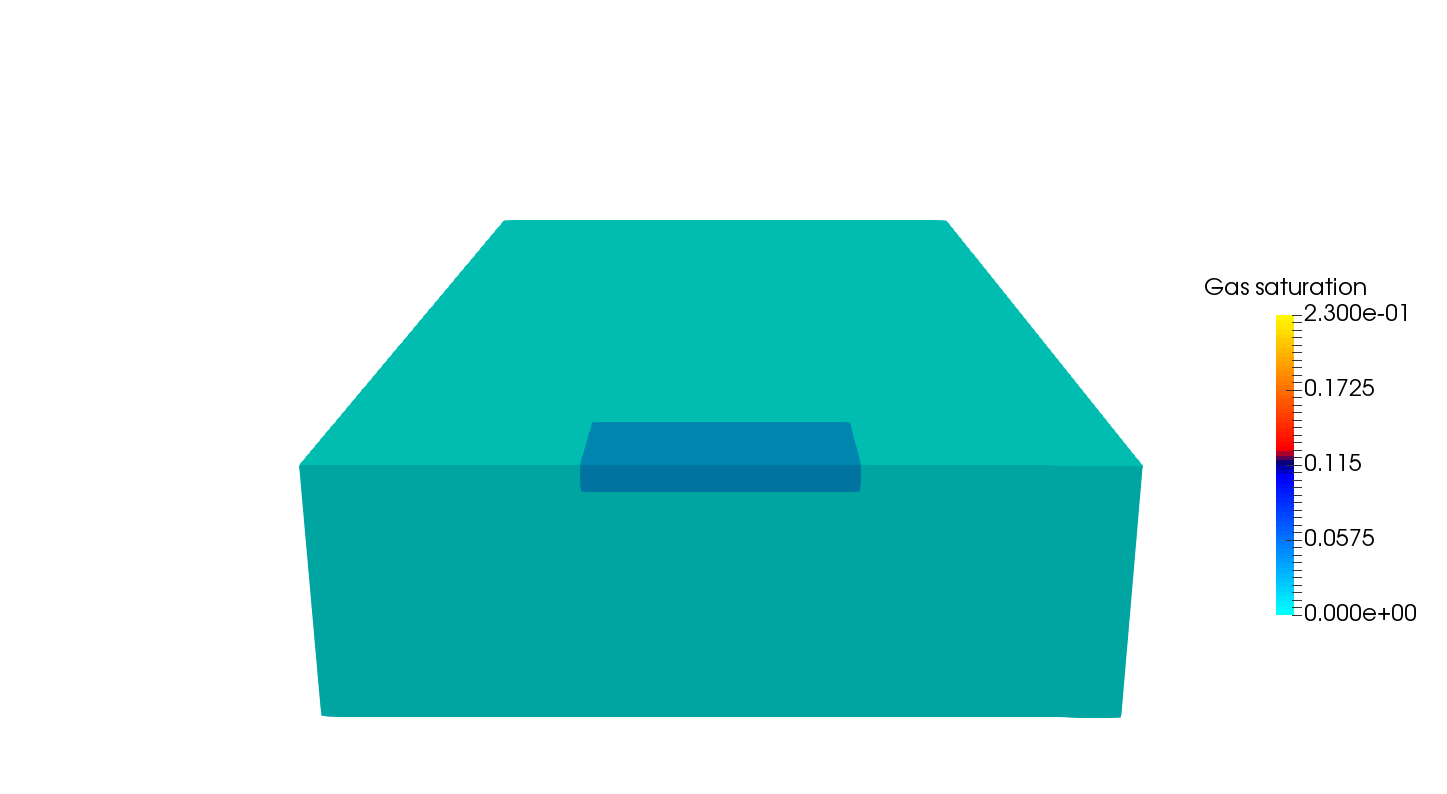}   
        \caption{Mesh size $h_1$.}
    \end{subfigure}
    \begin{subfigure}[b]{0.45\textwidth}
      \includegraphics[width=\textwidth]{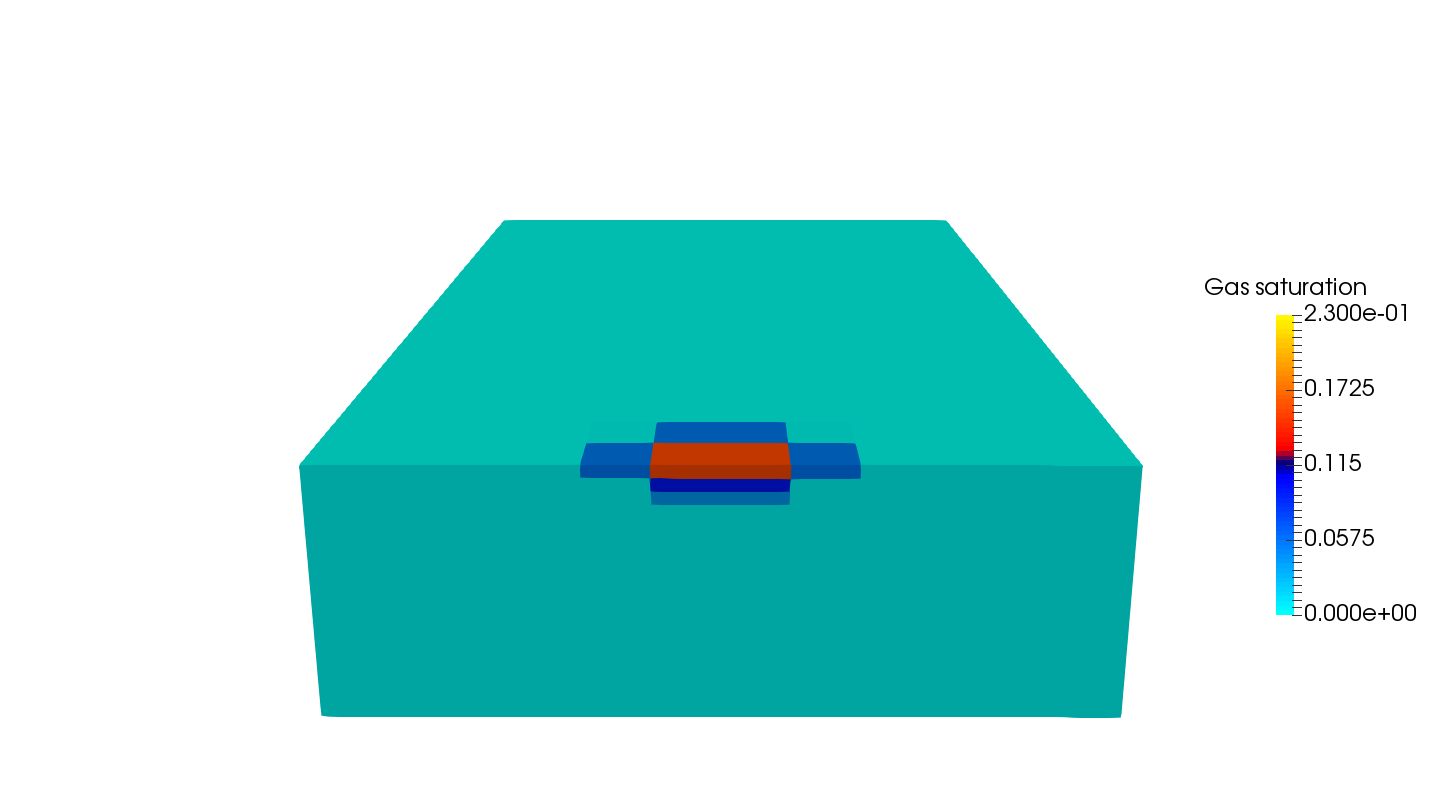}   
      \caption{Mesh size $h_2$.}
    \end{subfigure}\\
    \begin{subfigure}[b]{0.45\textwidth}
      \includegraphics[width=\textwidth]{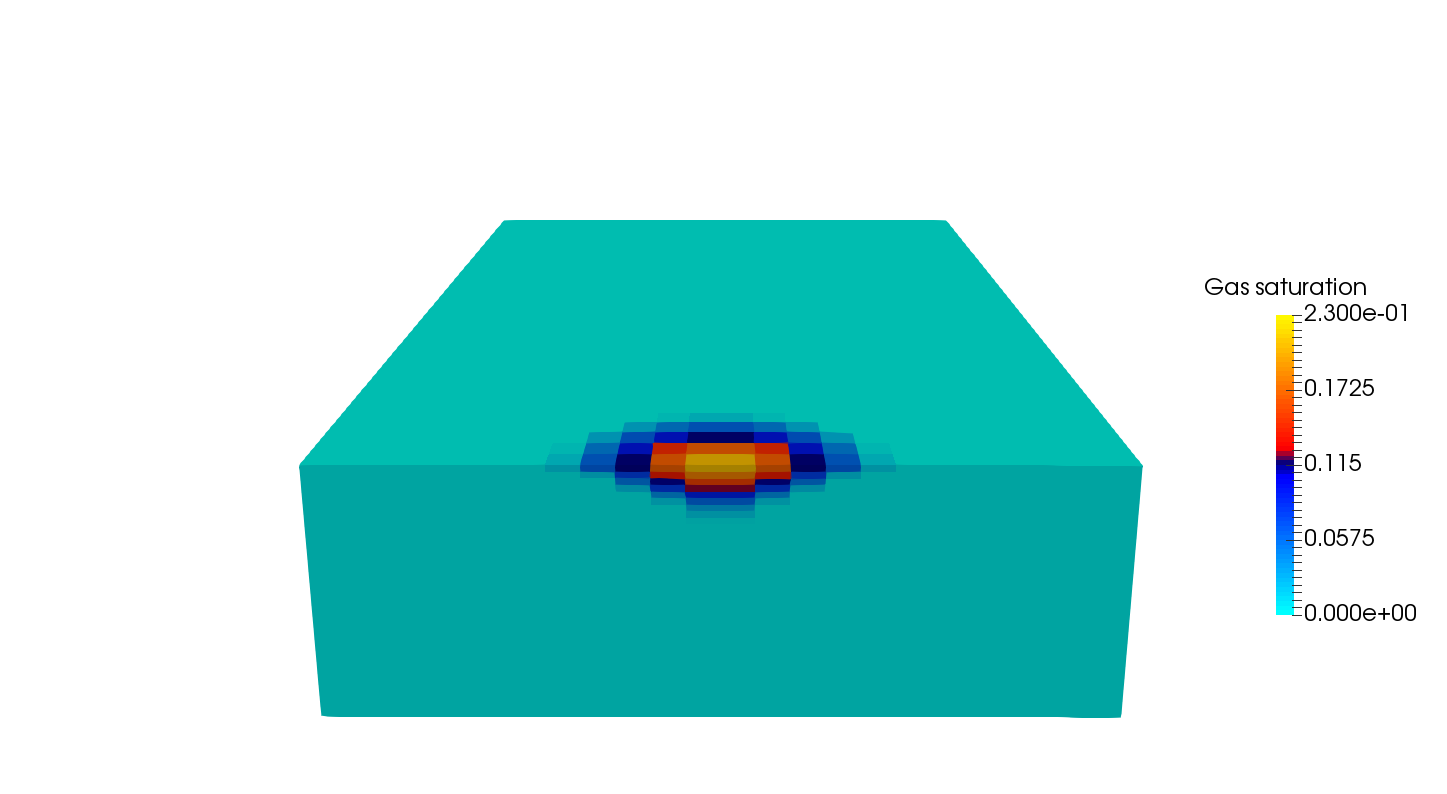}   
      \caption{Mesh size $h_3$.}
    \end{subfigure}
    \begin{subfigure}[b]{0.45\textwidth}
      \includegraphics[width=\textwidth]{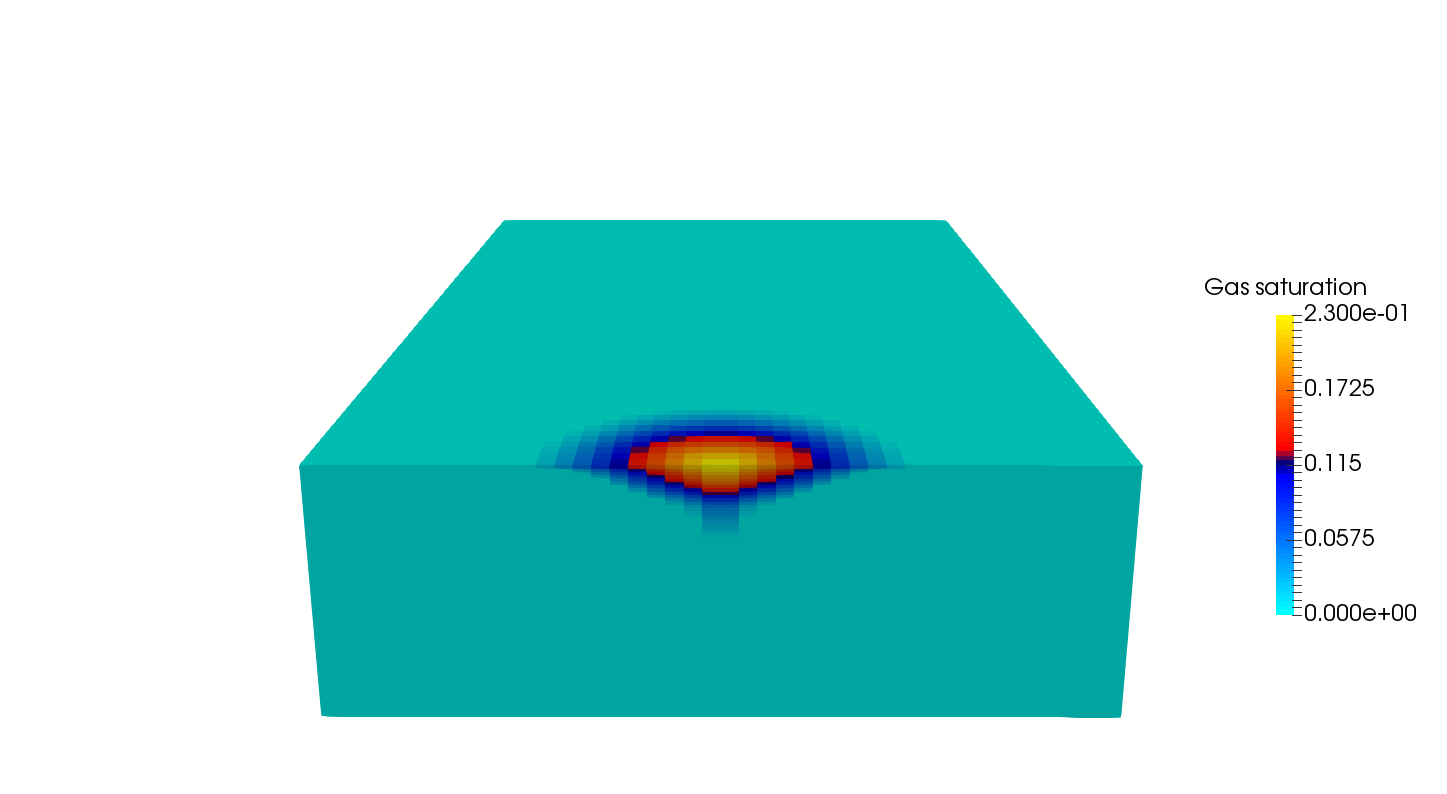}   
      \caption{Mesh size $h_4$.}
    \end{subfigure}
    \caption{Clip and close look of the gas saturation inside the reservoir at final time for all meshes (cell values).}
    \label{fig_gas_sat_res_all}
\end{figure}

Figures \ref{fig_gas_vol_well}
and \ref{fig_gas_vol_res}
show the total volume of gas
inside  the well,
and the total volume of gas
inside  the reservoir as functions of time 
for the family of refined meshes.
The solutions on the two coarsest meshes are still rough, which is expected given that the gas bubble is concentrated on a small region around the well (see Figure \ref{fig_gas_sat_res_all}). On the other hand the solutions on the two finest meshes are quite close exhibiting the good convergence of the scheme. 

\begin{figure}[H]
\begin{center}
\includegraphics[width=0.55\textwidth]{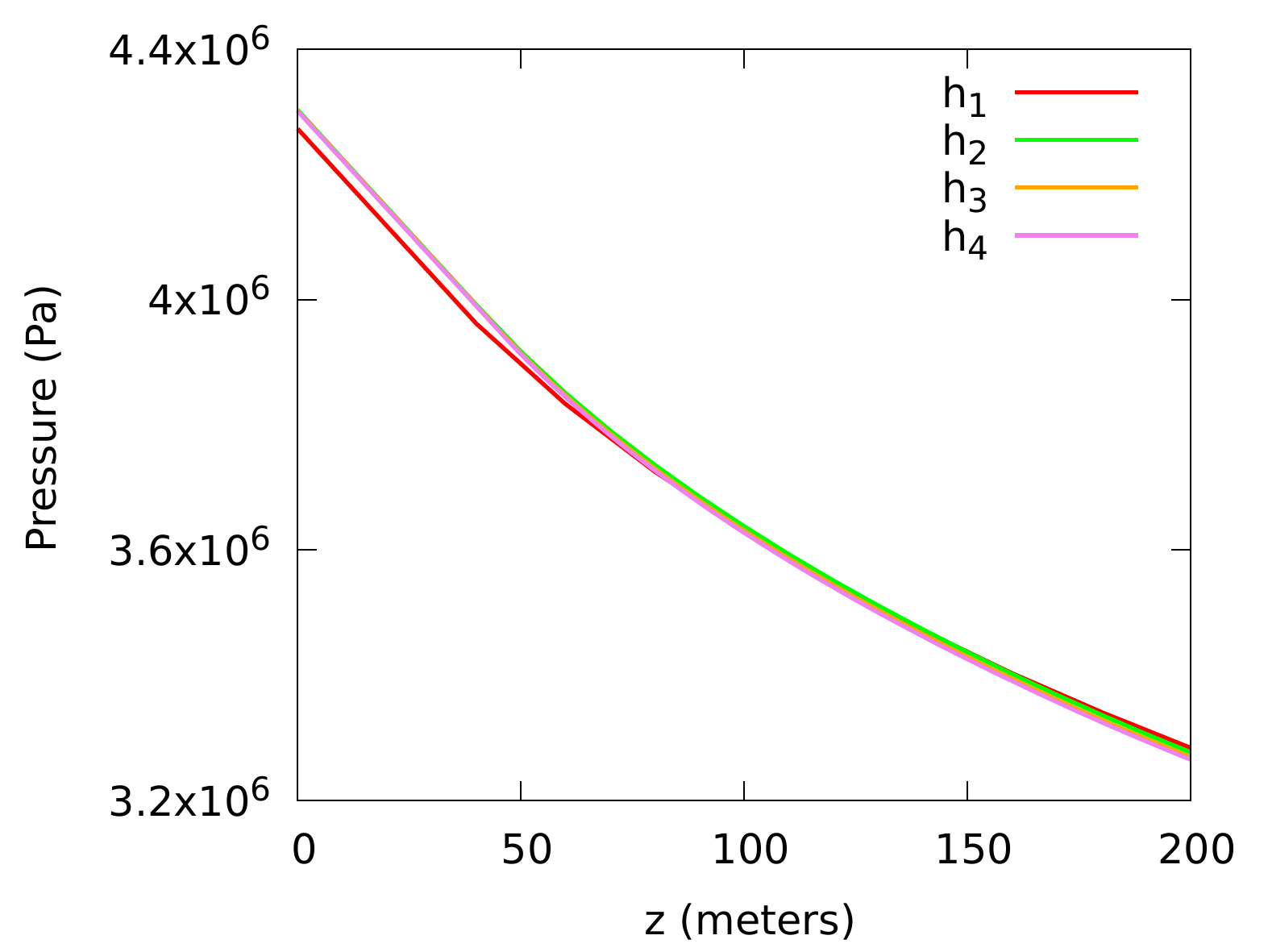}
    \caption{Pressure in Pa along the well at final time on the different meshes.}
    \label{fig_p_along_well}
 \end{center}
\end{figure}

\begin{figure}[H]
\begin{center}
\includegraphics[width=0.55\textwidth]{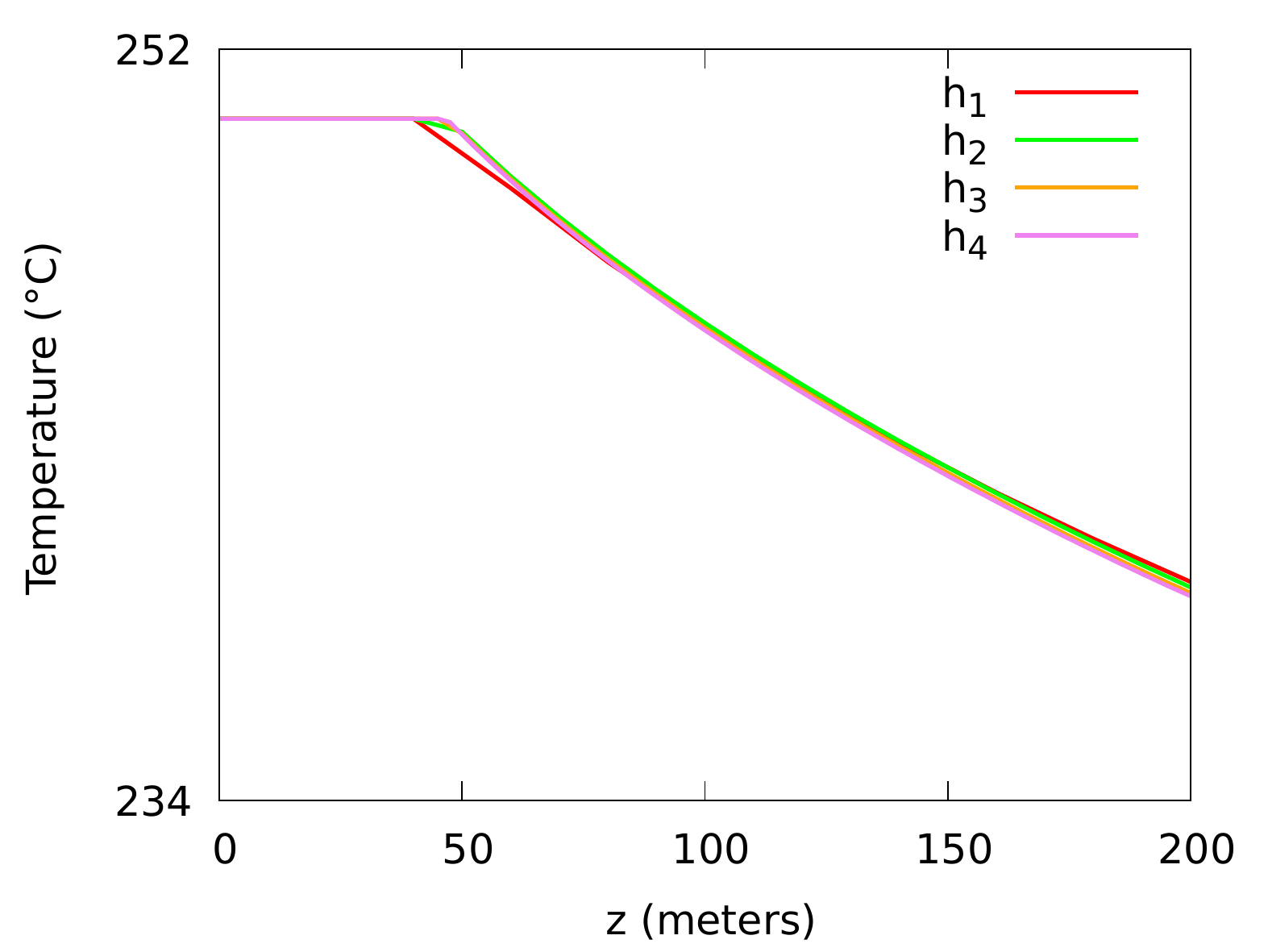}
    \caption{Temperature in $^\circ$C along the well at final time on the different meshes.}
    \label{fig_T_along_well}
  \end{center}
\end{figure}

\begin{figure}[H]
\begin{center}
\includegraphics[width=0.55\textwidth]{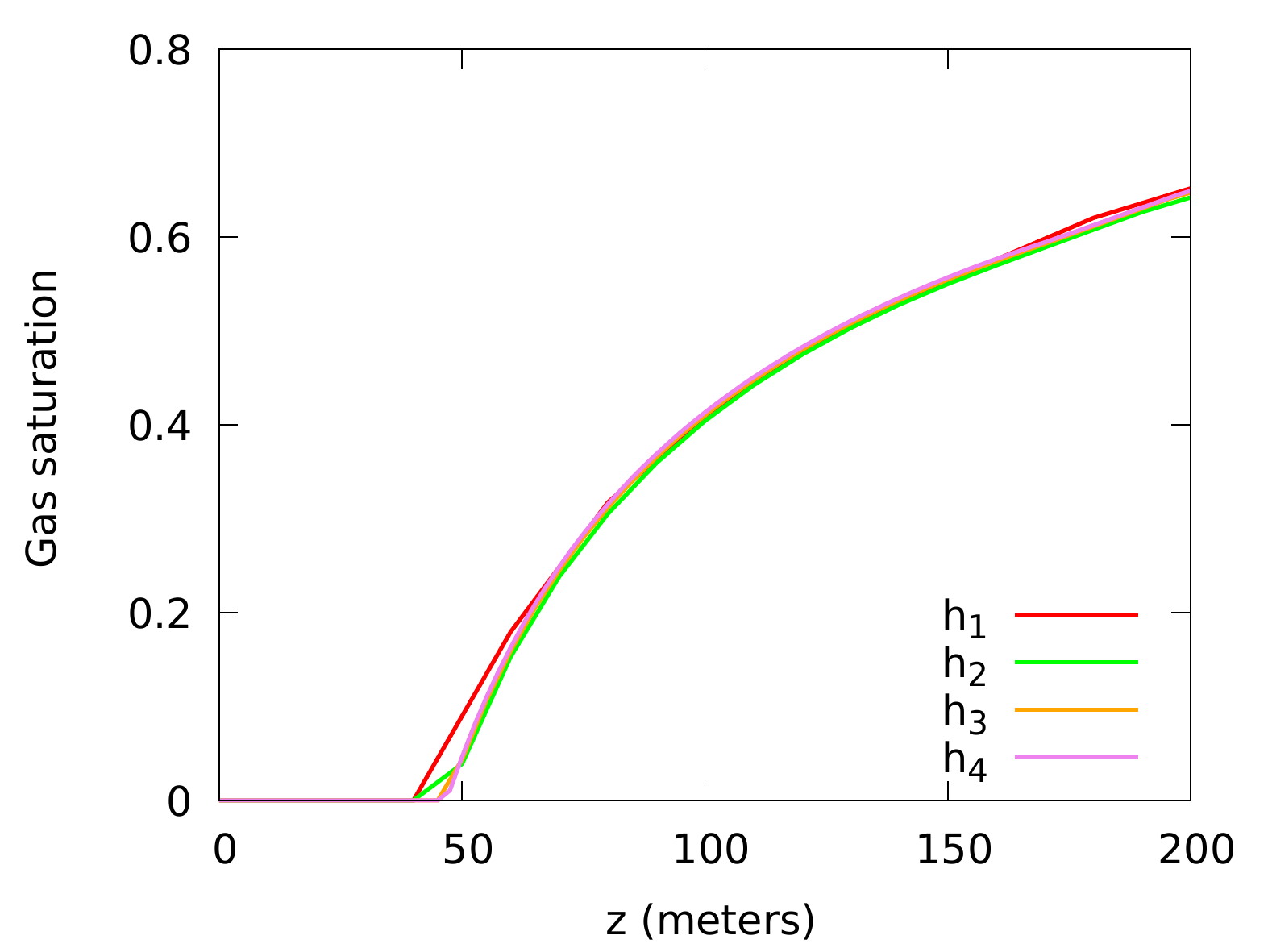}
     \caption{Gas saturation along the well  at final time on the different meshes.}
    \label{fig_g_sat_along_well}
 \end{center}
\end{figure}

In addition,
Figures \ref{fig_p_along_well}, \ref{fig_T_along_well},
and  \ref{fig_g_sat_along_well}
show the pressure, the temperature and the gas saturation along the well;
respectively, at final time $t_f$. The solutions are pretty close for all meshes and exhibit a good convergence behavior.

\begin{figure}[H]
    \centering
    \begin{subfigure}[b]{0.45\textwidth}
      \includegraphics[width=\textwidth]{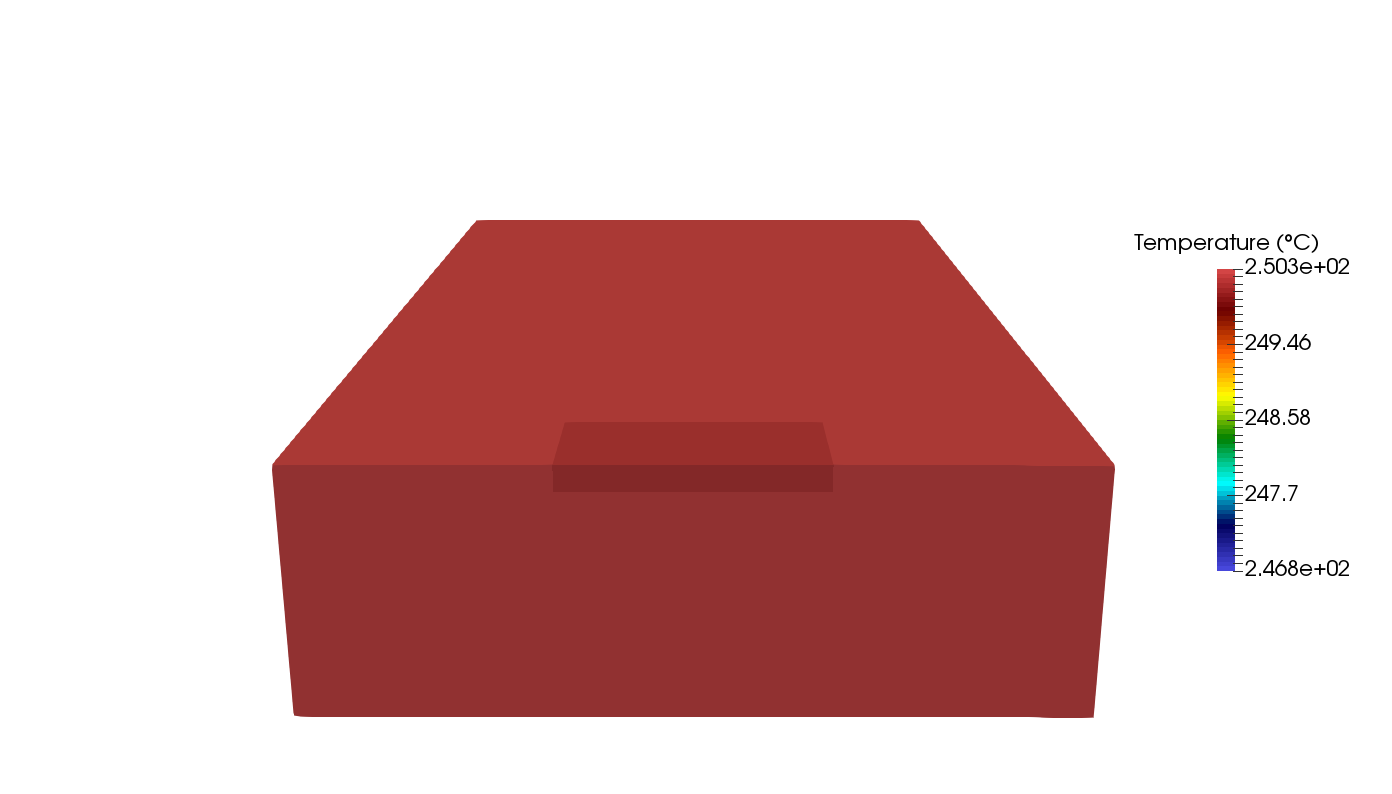}   
        \caption{Mesh size $h_1$.}
    \end{subfigure}
    \begin{subfigure}[b]{0.45\textwidth}
      \includegraphics[width=\textwidth]{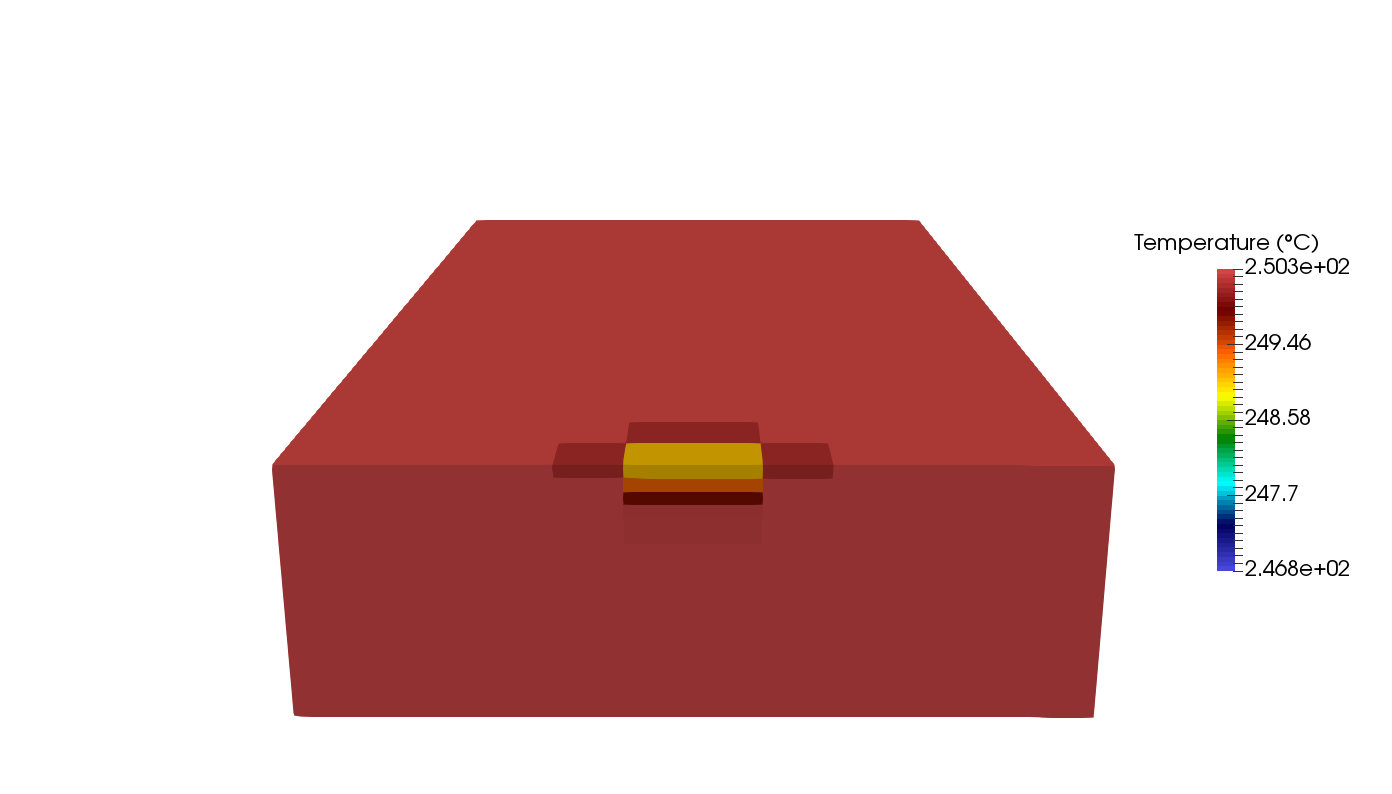}   
      \caption{Mesh size $h_2$.}
    \end{subfigure}\\
    \begin{subfigure}[b]{0.45\textwidth}
      \includegraphics[width=\textwidth]{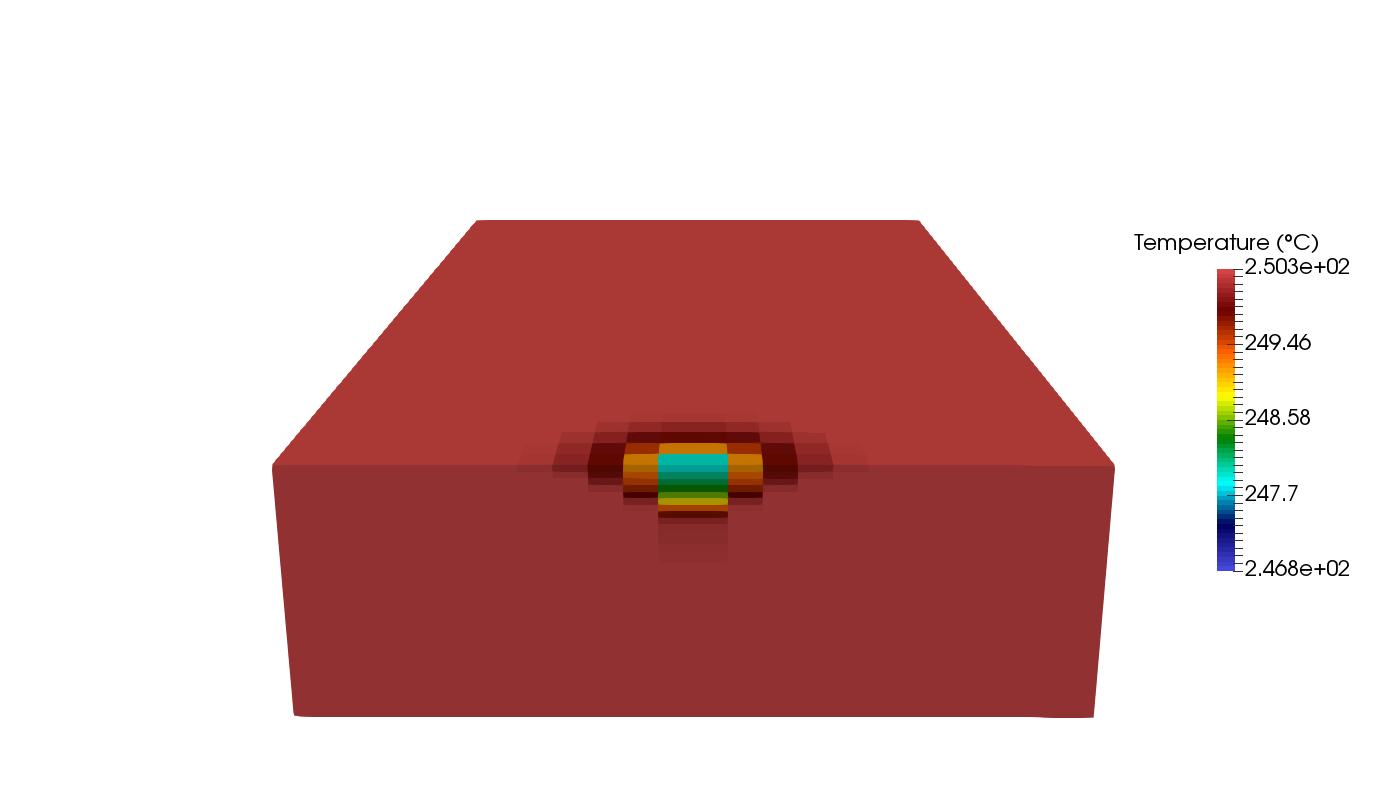}   
      \caption{Mesh size $h_3$.}
    \end{subfigure}
    \begin{subfigure}[b]{0.45\textwidth}
      \includegraphics[width=\textwidth]{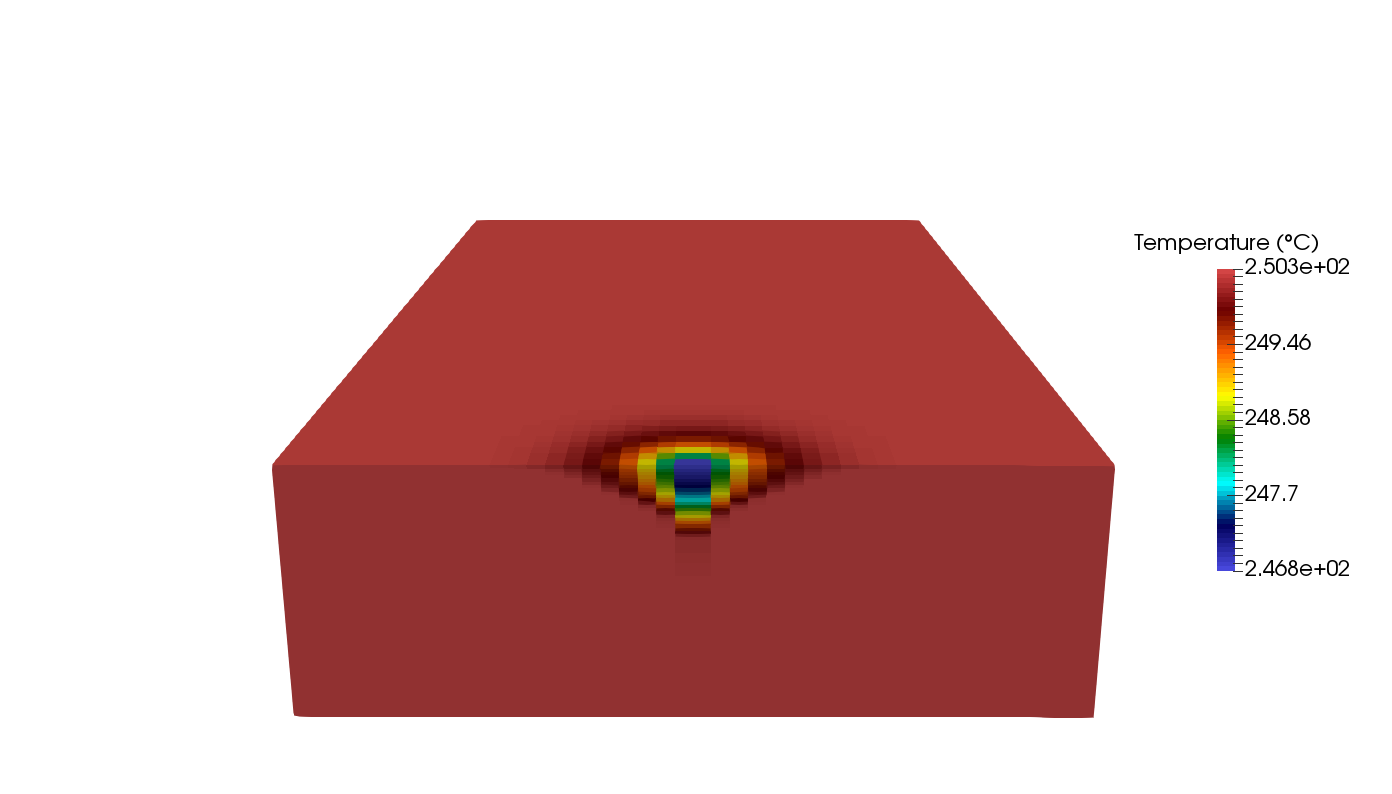}   
      \caption{Mesh size $h_4$.}
    \end{subfigure}
    \caption{Clip and close look of the temperature in $^\circ$C inside the reservoir at final time for all meshes (cell values).}
    \label{fig_T_res_all}
\end{figure}

%Figure \ref{fig_gas_vol_res}
%shows the gas saturation at final time $t_f$ inside the reservoir for the finest mesh $h_4$ exhibiting the cone shaped bubble of gas along the well at the top of the reservoir. 
Figures   \ref{fig_gas_sat_res_all},
and  \ref{fig_T_res_all}
show  a close look of the pressure and of the temperature inside the reservoir; respectively,
for all meshes at final time $t_f$. It illustrates the cone shaped bubble of gas along the well at the top of the reservoir and demontrates again the good convergence behavior of the discrete model.

At each time step, the nonlinear system is solved using a Newton algorithm.
The GMRES stopping criterion on the relative residual is fixed to $10^{-8}$.
The Newton solver is convergent if the relative residual is lower than $10^{-8}$ as well.

Table \ref{tab_numbeh} shows the numerical
efficiency of the proposed scheme for all meshes for the second stage of the simulation.
We denote by $N_{\Delta t}$ the number
of successful time steps,
by $N_{\text{Newton}}$ the average number of Newton iterations per
successful time step, and by $N_{\text{GMRES}}$ the average number of GMRES iterations per Newton iteration.
It exhibits a very good robustness of the Newton solver on the family of refined meshes and a moderate increase of the number of GMRES iterations with the mesh size.

\begin{figure}[H]
\begin{center}
\includegraphics[width=0.55\textwidth]{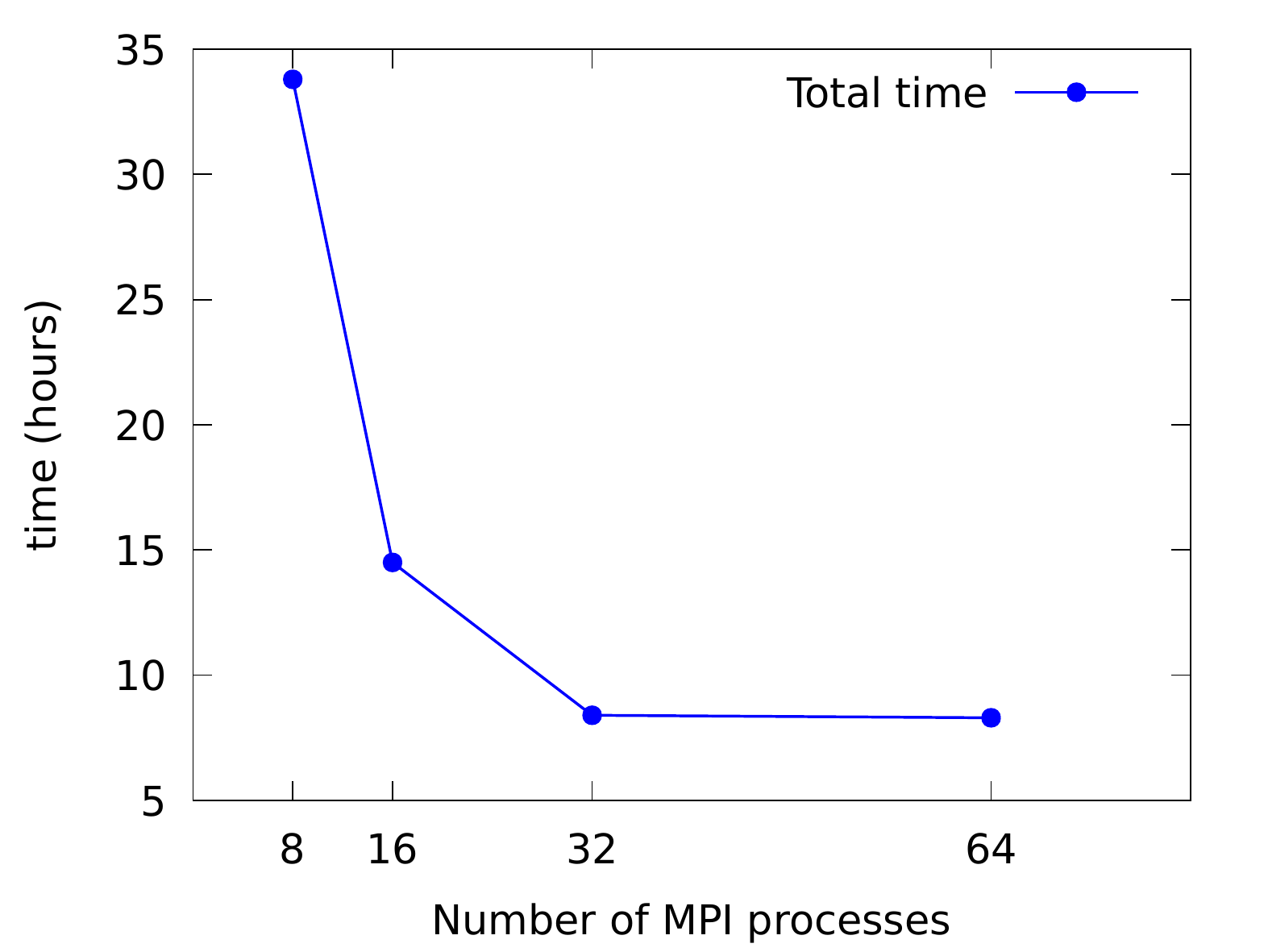}
    \caption{Total computational time vs. number of MPI processes for the second stage simulation on the finest mesh $h_4$.}
    \label{fig_speedup}
  \end{center}
\end{figure}

Finally, we present in Figure \ref{fig_speedup} the total computational time in hours obtained with the finest mesh $h_4$ for different numbers of MPI processes $N_p = 8, 16, 32, 64$. As usual for this type of simulations, the strong scalability is limited by the AMG preconditioner of the pressure block which requires a sufficiently high number of unknowns per processor to keep a good scalability, corresponding to roughly speaking $4~10^4$. This explains the good speed up obtained between $8$ and $32$ processors whereas the speed up becomes very small between $32$ and $64$ processors. 

\begin{table}[H]
\begin{center}
\begin{tabular}{|c|c|c|c|c|c|} \hline
  Mesh  & $\# \cells $  & $N_{\Delta t}$ &   $N_{\text{Newton}}$  &  $N_{\text{GMRES}}$  \\ \hline
  $h_1$ & 4000 & 134   & 1.99 & 8.59 \\
  $h_2$ & 32000 & 134   & 1.74 & 9.93 \\
  $h_3$ & 256000 & 134   & 1.92 & 11.75 \\
  $h_4$ & 1848320 & 133   & 2.22 & 15.91 \\\hline
\end{tabular}
\caption{Numerical behavior of the second stage of the simulation for different mesh sizes.
$N_{\Delta t}$ is the number
of successful time steps, $N_{\text{Newton}}$ the average number of Newton iterations per
successful time step, and $N_{\text{GMRES}}$ is the average number of GMRES iterations per Newton iteration.}  
  \label{tab_numbeh}
\end{center}
\end{table}

%\begin{figure}
 %   \centering
 %   \begin{subfigure}[b]{0.45\textwidth}
 %     
 %         \includegraphics[width=\textwidth]{ImagesCasTest/CasImage01.png} 
 %       \caption{Mesh size $h_1$.}
%   \end{subfigure}
%  \end{figure}
 
\subsection{Study of a high enthalpy reservoir}   
\begin{figure}[H]
    \centering
    \begin{subfigure}[b]{0.45\textwidth}
      \includegraphics[width=\textwidth]{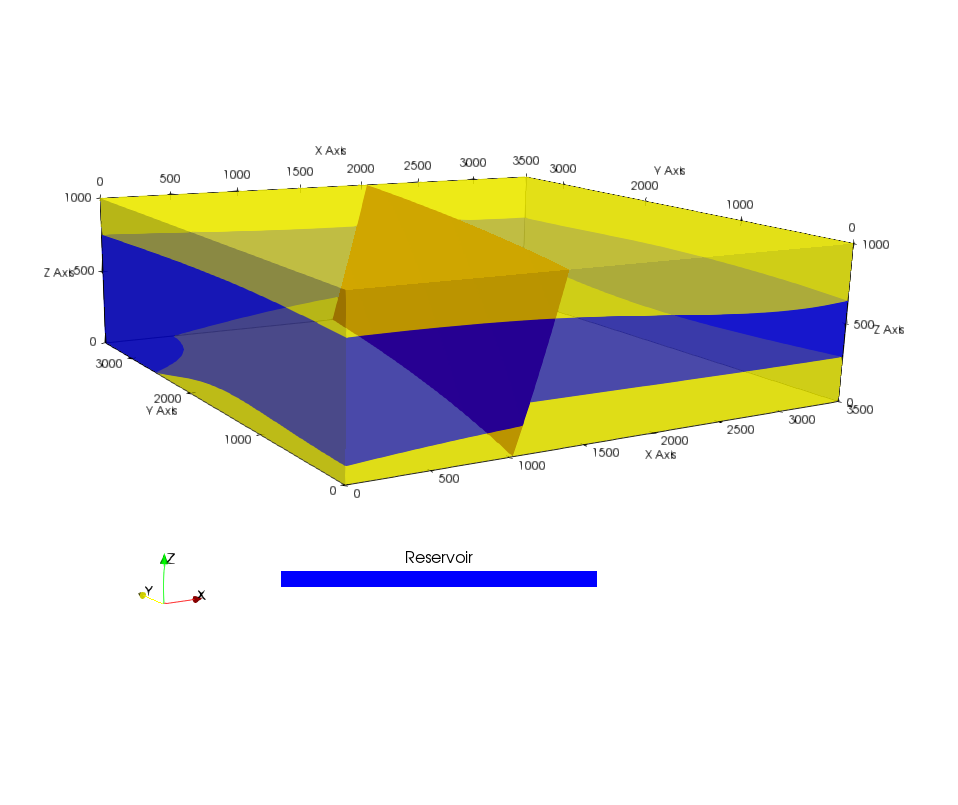}   
        \caption{Domain modelled.}
         \label{Domaincase2}
    \end{subfigure}
    \begin{subfigure}[b]{0.45\textwidth}
      \includegraphics[width=\textwidth]{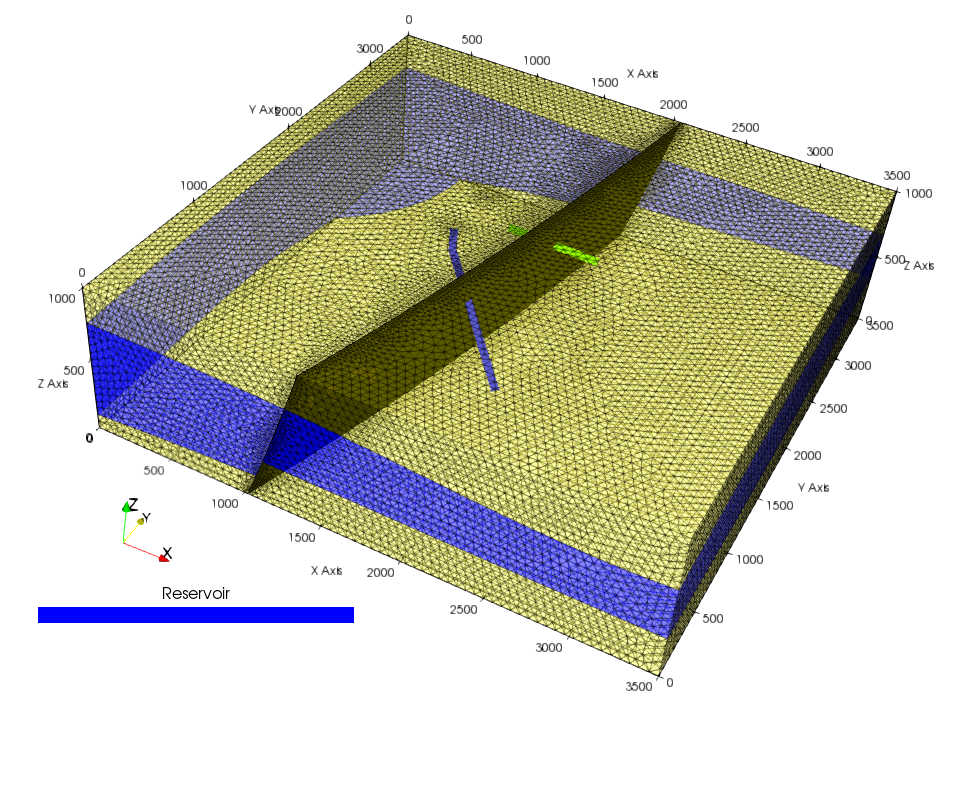}   
      \caption{Mesh and wells location.}
      \label{Meshcase2}
    \end{subfigure}
\caption{Geometry, mesh and wells data for the second numerical test.}
 \end{figure}

\hspace{\parindent} In this section, we consider a more realistic case built from geological and production data of a field in a volcanic area. %This case is inspired by the Bouillante geothermal field. The heat reservoir corresponds to a large rock volume affected by a normal fault. The extension of the heat reservoir is influenced by the geometry of the permeable fault.
The field is a convective dominated system initially in liquid phase, that is crossed by a major normal fault.\\
The reservoir (in blue in Figure \ref{Domaincase2}) is about $500$\;m thick; it is covered by  a weakly permeable clay caprock (in yellow) of $250$\;m thick, which outcrops at the surface. Below the reservoir is the basement layer (also in yellow).

% The rock mass is modelled up to the surface.
%Figure (\ref{Domaincase2}) presents the whole 3-D domain, which is a square box of 3500m side and 1000m high and decomposed into two kinds of matrices, one fault and two wellbores. The top of the reservoir (in blue, figure (\ref{Domaincase2})) is located at the base of the impermeable smectite caprock (in yellow,  (\ref{Domaincase2})) around 250m beneath the surface and is around 500m thick.\\
Figure \ref{Meshcase2} gives the tetrahedral mesh of the domain. The VAG finite volume discretization makes it possible to deal with complex geology including faults and complex well trajectories. The unstructured mesh of  700 000 tetrahedral elements draws on geological horizons. The fault is meshed as a two-dimensional (2D) surface, where the triangular elements are interconnected with the surrounding matrix using conformal meshing. The  (one-dimensional) wells  are discretized by a set of edges  as shown in Figure \ref{Meshcase2}. The computation of numerical well indexes would require an analytical solution for the linear diffusion equation, which is not known for
such a complex geometry involving fault and slanted wells. This solution could also be obtained numerically using a mesh at the scale of the wells, but its generation is out of the scope of this test case. Alternatively, we use for this test case an approximate analytical Peaceman type formula taking the fault into account and providing a good order of magnitude. \\

The geothermal field is operated using a doublet of two deviated wells, a producer (in green) and an injector (in blue), both of which cross the major fault as shown in Figure \ref{Meshcase2}.
The reservoir is assumed homogeneous with an isotropic permeability $\K_m = k_m I$, $k_m = 10^{-14}$\;m$^2$ and a porosity $\phi_m = 0.05$, while the faulted area has a thickness  $d_f = 10$ m, an isotropic permeability $\K_f = k_f I$, $k_f = 5.10^{-14}$\;m$^2$  and a porosity $\phi_f =0.05$.  The caprock and the basement layer are assumed weakly permeable  with $k_m = 10^{-19}$\;m$^2$.
The matrix and fracture thermal conductivities are set to $\lambda_m=\lambda_f=3$\;W.K$^{-1}$.m$^{-1}$ and 
the rock energy density is homogeneous for the whole rock mass such that $E_r(p,T) = \rho_{r} ~c^r_p ~T$ with $c^r_p= 1000$\;J.kg$^{-1}$.K$^{-1}$ and $\rho_{r}=2600$\;kg.m$^{-3}$.\\

\begin{figure}[H]
\centering
\includegraphics[width=13.5cm]{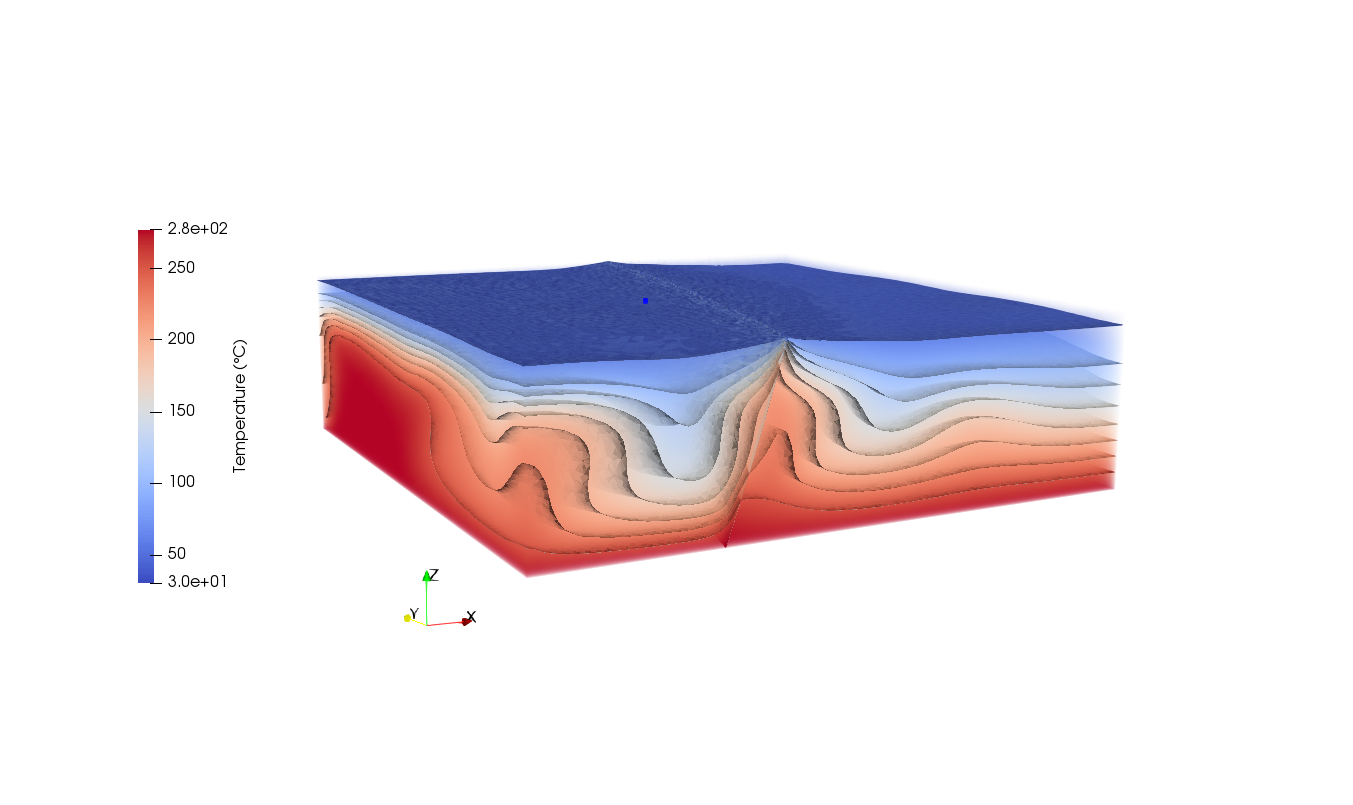}
\caption{Initial state dominated by convention. Isotemperature surfaces.}
\label{cas2_initial}
\end{figure}
As the previous numerical test, this simulation consists in two stages. The first one acts as a preliminary step where the initial state of the geothermal system, which is already dynamic, is achieved by performing a simulation over a long period (here $10^5$\;years) from an hydrostatic pressure state (with $1$\;bar at the top of the model), and a temperature field increasing linearly with depth (between $30\;^\circ$C at the top to $280\,^{\circ}$C at the bottom). Dirichlet boundary conditions for temperature and pressure are thus imposed at the top and bottom boundaries. No flow  and Dirichlet temperature conditions are applied on the lateral boundaries.
 The initial state obtained is convective; the fluid in the reservoir is in liquid state with a low fraction of gas near the top of the reservoir. Iso-temperature contours are represented in Figure \ref{cas2_initial} and show the development of convection cells and the influence of the fault, which is a more permeable zone.

\begin{figure}[H]
    \begin{center}
      \includegraphics[width=13cm]{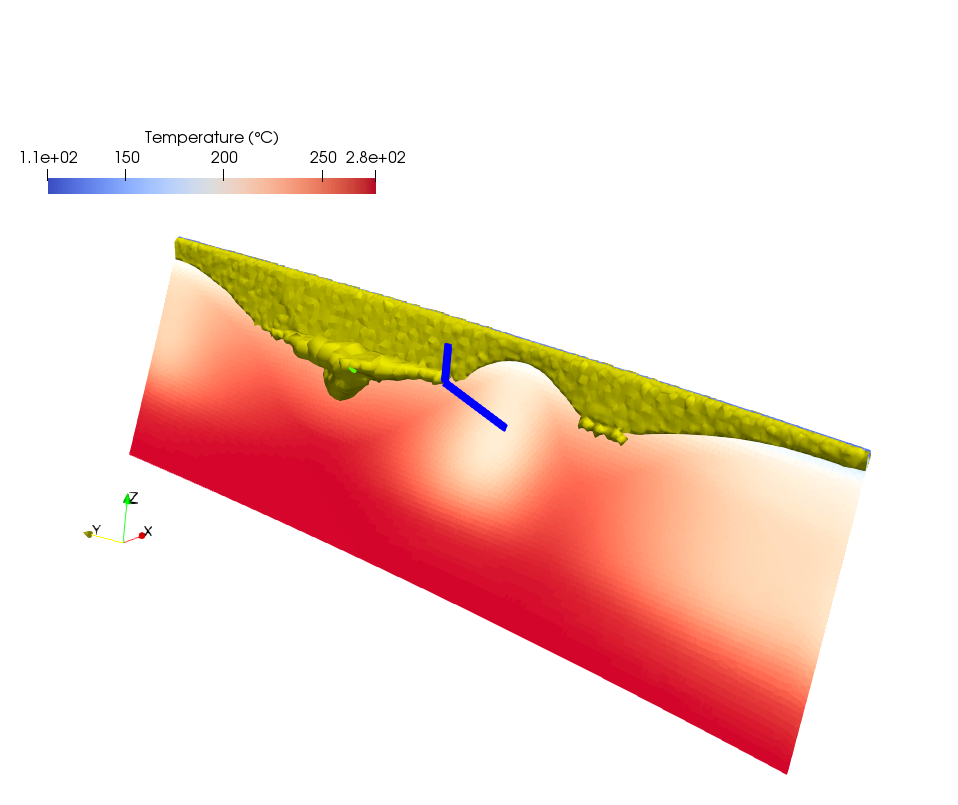}
    \caption{Temperature and saturation after $5$\;years of production - cells with a gas saturation greater than $0.1$ are filled in yellow - the temperature is represented in the fault plane}
    \label{cas2_afterproduction}
    \end{center}
\end{figure}

Then the second stage begins where the reservoir production starts with steam production at the producer well-head: 
a flow rate of $250$\;ton.hr$^{-1}$ is imposed at the well-head for five years. The same boundary conditions are imposed as in the initial state determination, but the temperature imposed on the lateral boundaries is now given by the average temperature distribution in the rock mass at this initial state.
The depletion occuring near the producer well favors the development of a steam cap in the reservoir as well as in the fault zone. Figure \ref{cas2_afterproduction} shows this steam cap: faces in the fault and cells in the reservoir with a gas saturation greater than 0.1 are filled in yellow, while temperature field is also represented on the other faces of the fault plane.
%The depletion occuring near the producer well favors the development of a steam cap in the reservoir as well as in the fault zone as shown in Figure (\ref{cas2_afterproduction}) where cells with a gas saturation greater than 0.1 are filled in yellow, while temperature field is represented on the fault plane. 

\begin{figure}[H]
    \begin{center}
      \includegraphics[width=13cm]{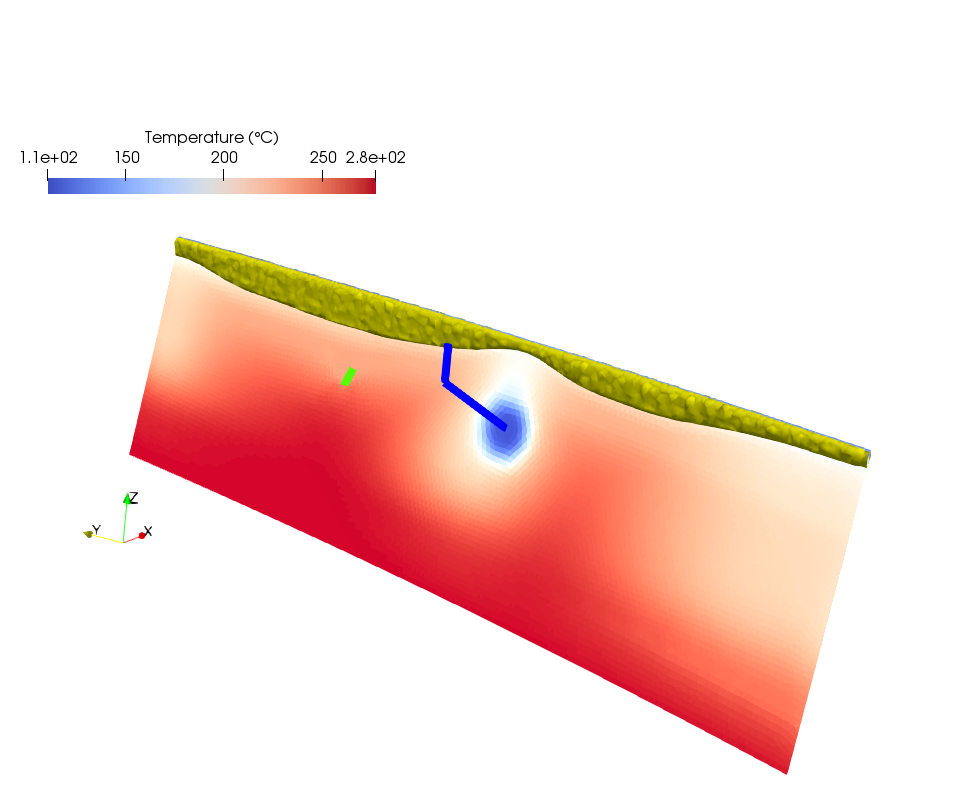}
    \caption{Temperature and saturation after $10$ years of production with reinjection during the last $5$ years  - cells with a gas saturation greater than $0.1$ are filled in yellow - the temperature is represented in the fault plane}
    \label{cas2_afterinjection}
    \end{center}
\end{figure}

\begin{figure}[H]
    \centering
    \begin{subfigure}[b]{0.45\textwidth}
      \includegraphics[width=\textwidth]{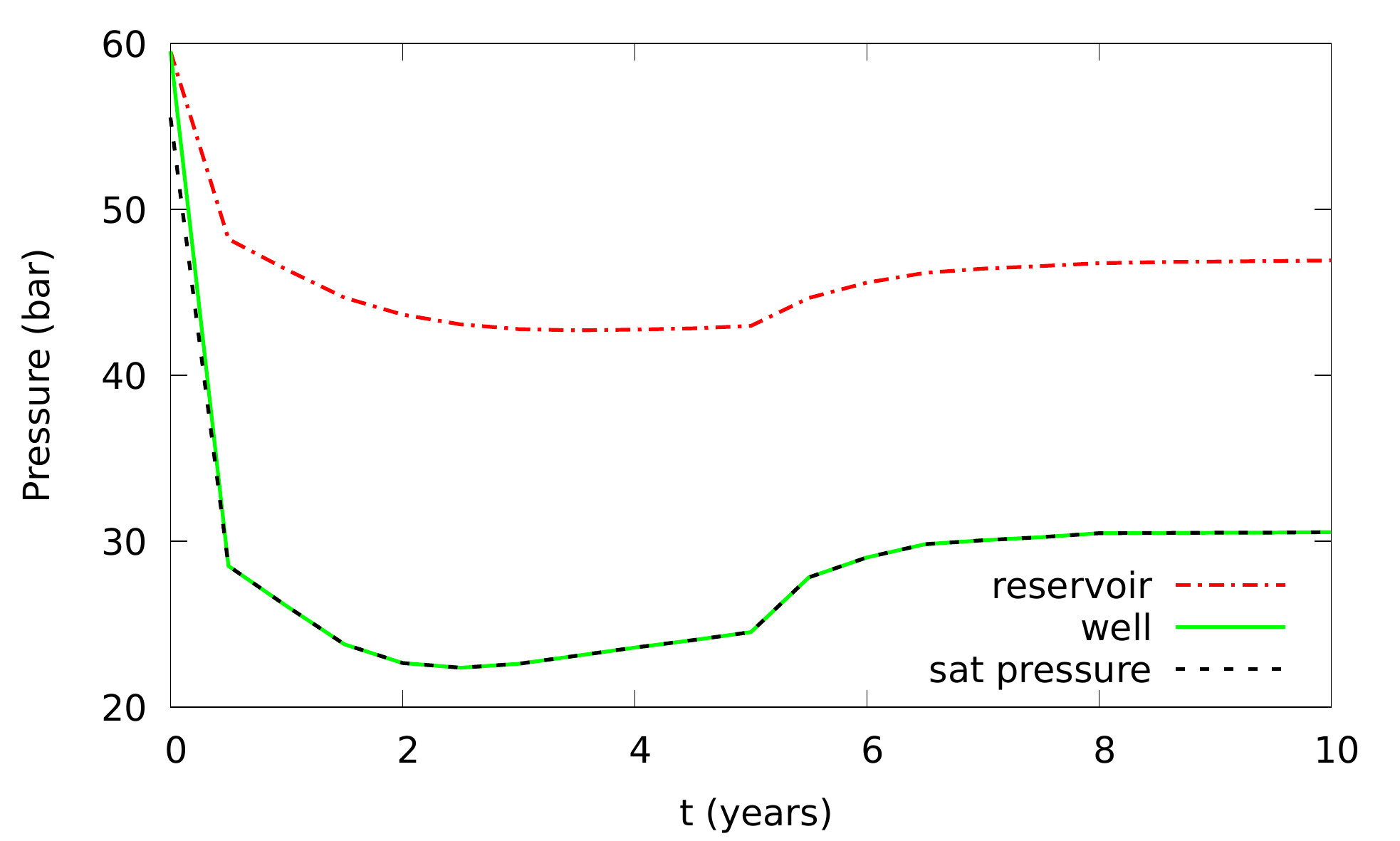}   
        \caption{Pressure evolution.}
         \label{WellPressure}
    \end{subfigure}
    \begin{subfigure}[b]{0.45\textwidth}
      \includegraphics[width=\textwidth]{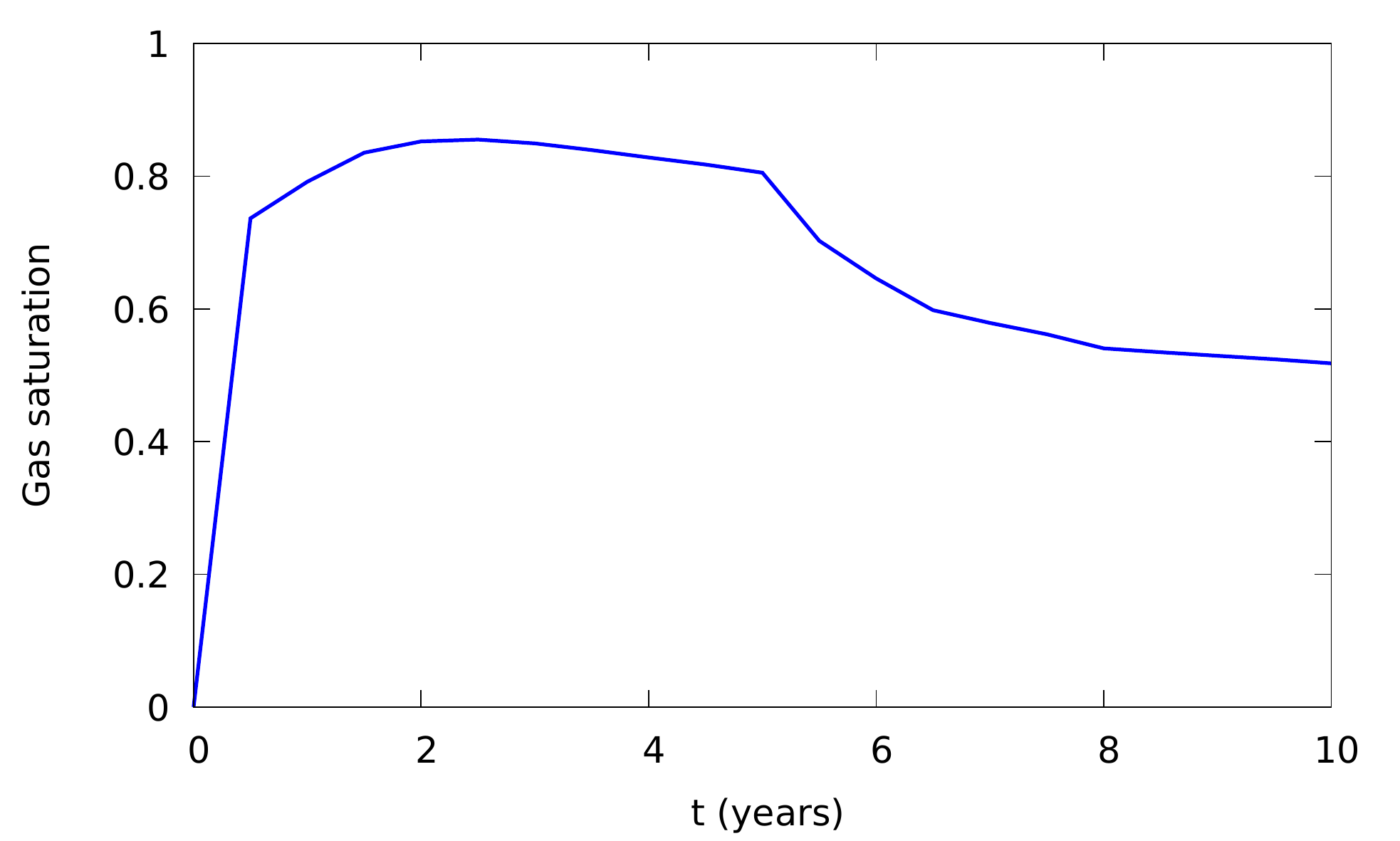}   
      \caption{Saturation evolution.}
      \label{WellSaturation}
    \end{subfigure}
     \caption{(a) Pressure evolution in the reservoir (red dashed line) and in the well (green line) at $455$\;m depth. Saturation pressure in the well at $455$\;m depth is given by the black dotted line. (b): Saturation evolution in the well at $455$\;m depth.}
 \end{figure}

 After five years of production and reservoir depletion, half of the fluid produced is reinjected at the injector with a wellhead temperature of $110\,^\circ$C. During the injection, vapor around the injector condenses and the steam cap generated around the producer is considerably reduced (Figure \ref{cas2_afterinjection}). 
%\begin{figure}[htbp!]
%    \begin{center}
 %     \def\svgscale{0.5}
 %    \includegraphics[width=\textwidth]{{ImagesCasTest/CasImage05.png}
 %   \caption{Initial state dominated by convention.}
 %   \label{cas2_maillage}
%    \end{center}
%\end{figure}
Figures \ref{WellPressure} and \ref{WellSaturation} show at a given depth of $455$\;m respectively the evolution of pressure in the reservoir and in the well and the saturation evolution in the well. Reservoir pressure decreases during the first five years of production, while reinjection of half of the fluid produced during the next five years leads to a pressure build-up in the reservoir (the model is not hydraulically closed). Well pressure follows the same trends. Whereas gas saturation was around 80$\%$  during the depletion phase in the well at $455$\;m depth, injection results in a reduced gas saturation in the well down to say 50$\%$ at $455$\;m depth (Figure \ref{WellSaturation}).
Table \ref{TableSpeedup} shows the numerical
efficiency of the proposed scheme  for both  stages of the simulation
and different numbers of MPI processes $N_p = 4, 8,16$. 
We use the same notations as in the previous test case and report in addition the total simulation time in hours. These results exhibit the very good robustness of the linear and nonlinear solvers w.r.t. the number of MPI processes. A very good speedup is obtained up to 16 MPI processes verifying that
parallel computing  makes  possible to have reasonable computation times to model industrial cases such as the one presented in this section. 

\begin{table}[H]
\begin{center}
\begin{tabular}{|c|c|c|c|c|c|} \hline
  Stage & $N_p$  &   $N_{\Delta t}$  &  $N_{\text{Newton}}$  &  $N_{\text{GMRES}}$ & Time (hrs)  \\ \hline
 \multirow{3}{*}{1} & 4 &  1515 & 4.6 & 29.3 & 98.2 \\
  & 8 &  1507 & 4.6 & 29.4 & 31.9 \\
  & 16 & 1526 & 4.6 & 30.0 & 17.8 \\\hline
 \multirow{3}{*}{2} &  4 &  1395 & 7.3 & 7.7  & 65.9 \\
  & 8 &  1367 & 7.2 & 7.6  & 20.2 \\
  & 16 & 1320 & 7.2 & 7.9  & 10.1 \\\hline
\end{tabular}
  \caption{Numerical behavior of both stages of the simulation for different number of processors. $N_{\Delta t}$ is the number
of successful time steps, $N_{\text{Newton}}$ the average number of Newton iterations per
successful time step, $N_{\text{GMRES}}$ the average number of GMRES iterations per Newton iteration, and Time (hrs) is the total simulation time in hours.}  
\label{TableSpeedup}
\end{center}
\end{table}

\section{Conclusion}
\label{sec_conclu}
This paper focuses on the numerical modelling of geothermal systems in complex geological settings.
The proposed approach is based on unstructured meshes to model complex features such as faults and deviated wells. It solves liquid vapor two-phase Darcy flows coupled with energy transfers and thermodynamical equilibrium. The use of the hybrid-dimensional polytopal VAG scheme allows to treat physically complex cases, while respecting geometrical constraints.  We particularly focus on the well modelling with deviated or multi-branch wells defined as a collection of edges of the mesh with rooted tree data structure. By using an explicit pressure drop calculation, the well model reduces to a single equation with only one well implicit unknown fully coupled to the reservoir system. Finally, efficient parallel linear and nonlinear solvers ensure acceptable computation times on real case studies.
A sanity checked is first presented showing the numerical convergence of the discrete model  on a diphasic vertical producer well in a simple reservoir geometry. Then, the efficiency of our approach is demonstrated on a geothermal test case of high enthalpy faulted reservoir using a doublet of two deviated wells crossing the fault.

An improved model of cross flows between well and reservoir will be investigated in the near future. Industrial studies of high and medium enthalpy geothermal reservoirs are currently under way with the approach proposed in this paper.

\section*{Acknowledgments} 

\thanks{This work was supported by a joint project between Storengy, BRGM and UCA and by the CHARMS ANR project (ANR-16-CE06-0009).}

%%-----------------------------
%%      your bibliography
%%-----------------------------

\bibliographystyle{abbrv}
%\bibliography{Bib-rev}
 \bibliography{article}

\begin{thebibliography}{10}

\bibitem{Aav03}
I.~Aavatsmark and R.~Klausen.
\newblock {Well Index in Reservoir Simulation for Slanted and Slightly Curved
  Wells in 3D Grids}.
\newblock {\em SPE Journal}, 8(01):41--48, 03 2003.

\bibitem{AELHP152D}
R.~Ahmed, M.~Edwards, S.~Lamine, B.~Huisman, and M.~Pal.
\newblock {Control-volume distributed multi-point flux approximation coupled
  with a lower-dimensional fracture model}.
\newblock {\em Journal of Computational Physics}, 284:462--489, mar 2015.

\bibitem{AELHP153D}
R.~Ahmed, M.~G. Edwards, S.~Lamine, B.~A. Huisman, and M.~Pal.
\newblock {Three-dimensional control-volume distributed multi-point flux
  approximation coupled with a lower-dimensional surface fracture model}.
\newblock {\em Journal of Computational Physics}, 303:470--497, dec 2015.

\bibitem{MAE02}
C.~Alboin, J.~Jaffr{\'{e}}, J.~Roberts, and C.~Serres.
\newblock {Modeling fractures as interfaces for flow and transport in porous
  media}.
\newblock volume 295, pages 13--24, 2002.

\bibitem{ABH09}
P.~Angot, F.~Boyer, and F.~Hubert.
\newblock {Asymptotic and numerical modelling of flows in fractured porous
  media}.
\newblock {\em ESAIM: Mathematical Modelling and Numerical Analysis},
  43(2):239--275, mar 2009.

\bibitem{AFSVV16}
P.~F. Antonietti, L.~Formaggia, A.~Scotti, M.~Verani, and N.~Verzott.
\newblock {Mimetic finite difference approximation of flows in fractured porous
  media}.
\newblock {\em ESAIM M2AN}, 50:809--832, 2016.

\bibitem{Aunzo91}
Z.~P. Aunzo, G.~Bjornsson, and G.~S. Bodvarsson.
\newblock {Wellbore Models GWELL, GWNACL, and HOLA, user's guide}.
\newblock Technical Report LBL-31428, Earth Sciences Division, Lawrence
  Berkeley National Laboratory, University of California, 1991.

\bibitem{aziz-settari-79}
K.~Aziz and A.~Settari.
\newblock {\em {Petroleum Reservoir Simulation}}.
\newblock Applied Science Publishers, 1979.

\bibitem{Xing2018}
{Beaude, Laurence}, {Beltzung, Thibaud}, {Brenner, Konstantin}, {Lopez, Simon},
  {Masson, Roland}, {Smai, Farid}, {Thebault, Jean-fr\'ed\'eric}, and {Xing,
  Feng}.
\newblock Parallel geothermal numerical model with fractures and multi-branch
  wells.
\newblock {\em ESAIM: ProcS}, 63:109--134, 2018.

\bibitem{BMTA03}
I.~I. Bogdanov, V.~V. Mourzenko, J.-F. Thovert, and P.~M. Adler.
\newblock {Two-phase flow through fractured porous media}.
\newblock {\em Physical Review E}, 68(2), aug 2003.

\bibitem{GSDFN}
K.~Brenner, M.~Groza, C.~Guichard, G.~Lebeau, and R.~Masson.
\newblock {Gradient discretization of hybrid-dimensional Darcy flows in
  fractured porous media}.
\newblock {\em Numerische Mathematik}, 134(3):569--609, nov 2016.

\bibitem{BGGM14}
K.~Brenner, M.~Groza, C.~Guichard, and R.~Masson.
\newblock {Vertex Approximate Gradient scheme for hybrid-dimensional two-phase
  Darcy flows in fractured porous media}.
\newblock {\em ESAIM: Mathematical Modelling and Numerical Analysis},
  2(49):303--330, 2015.

\bibitem{BGJMP17}
K.~Brenner, M.~Groza, L.~Jeannin, R.~Masson, and J.~Pellerin.
\newblock Immiscible two-phase {D}arcy flow model accounting for vanishing and
  discontinuous capillary pressures: application to the flow in fractured
  porous media.
\newblock {\em Computational Geosciences}, 21(5):1075--1094, 2017.

\bibitem{Brenner2021}
K.~Brenner, J.~Hennicker, and R.~Masson.
\newblock {\em Nodal Discretization of Two-Phase Discrete Fracture Matrix
  Models}, pages 73--118.
\newblock Springer International Publishing, Cham, 2021.

\bibitem{BHMS2016}
K.~Brenner, J.~Hennicker, R.~Masson, and P.~Samier.
\newblock {Gradient discretization of hybrid-dimensional Darcy flow in
  fractured porous media with discontinuous pressures at matrix-fracture
  interfaces}.
\newblock {\em IMA Journal of Numerical Analysis}, sep 2016.

\bibitem{BHMS18}
K.~Brenner, J.~Hennicker, R.~Masson, and P.~Samier.
\newblock Hybrid-dimensional modelling of two-phase flow through fractured
  porous media with enhanced matrix fracture transmission conditions.
\newblock {\em Journal of Computational Physics}, 357:100--124, 2018.

\bibitem{CZ09}
Z.~Chen and Y.~Zhang.
\newblock {Well flow models for various numerical methods}.
\newblock {\em J. Numer. Anal. Model.}, 6:375--388, 2009.

\bibitem{Eymard.Herbin.ea:2010}
R.~Eymard, C.~Guichard, and R.~Herbin.
\newblock {Small-stencil 3D schemes for diffusive flows in porous media}.
\newblock {\em ESAIM: Mathematical Modelling and Numerical Analysis},
  46(2):265--290, 2012.

\bibitem{EHGM-CG-12}
R.~Eymard, C.~Guichard, R.~Herbin, and R.~Masson.
\newblock {Vertex-centred discretization of multiphase compositional Darcy
  flows on general meshes}.
\newblock {\em Computational Geosciences}, 16(4):987--1005, 2012.

\bibitem{FFJR16}
I.~Faille, A.~Fumagalli, J.~Jaffr{\'{e}}, and J.~E. Roberts.
\newblock {Model reduction and discretization using hybrid finite volumes of
  flow in porous media containing faults}.
\newblock {\em Computational Geosciences}, 20:317--339, 2016.

\bibitem{FNFM03}
E.~Flauraud, F.~Nataf, I.~Faille, and R.~Masson.
\newblock {Domain decomposition for an asymptotic geological fault modeling}.
\newblock {\em Comptes Rendus M{\'{e}}canique}, 331(12):849--855, dec 2003.

\bibitem{Gjerde20}
I.~G. Gjerde, K.~Kumar, and J.~M. Nordbotten.
\newblock A singularity removal method for coupled 1d--3d flow models.
\newblock {\em Computational Geosciences}, 24(2):443--457, 2020.

\bibitem{GRANET200135}
S.~Granet, P.~Fabrie, P.~Lemonnier, and M.~Quintard.
\newblock A two-phase flow simulation of a fractured reservoir using a new
  fissure element method.
\newblock {\em Journal of Petroleum Science and Engineering}, 32(1):35 -- 52,
  2001.

\bibitem{HADEH09}
H.~Haegland, A.~Assteerawatt, H.~Dahle, G.~Eigestad, and R.~Helmig.
\newblock {Comparison of cell- and vertex-centered discretization methods for
  flow in a two-dimensional discrete-fracture-matrix system}.
\newblock {\em Advances in Water resources}, 32:1740--1755, 2009.

\bibitem{HY:2001}
V.~E. Henson and U.~M. Yang.
\newblock {BoomerAMG: A parallel algebraic multigrid solver and
  preconditioner}.
\newblock {\em Applied Numerical Mathematics}, 41(1):155--177, 2002.

\bibitem{HF08}
H.~Hoteit and A.~Firoozabadi.
\newblock {An efficient numerical model for incompressible two-phase flow in
  fractured media}.
\newblock {\em Advances in Water Resources}, 31(6):891--905, jun 2008.

\bibitem{KDA04}
M.~Karimi-Fard, L.~Durlofsky, and K.~Aziz.
\newblock {An efficient discrete-fracture model applicable for general-purpose
  reservoir simulators}.
\newblock {\em SPE Journal}, 9(02):227--236, jun 2004.

\bibitem{Kr11}
S.~Kr\"autle.
\newblock The semi-smooth newton method for multicomponent reactive transport
  with minerals.
\newblock {\em Advances in Water Resources}, 34:137--151, 2011.

\bibitem{LVW:2001}
S.~Lacroix, Y.~V. Vassilevski, and M.~F. Wheeler.
\newblock {Decoupling preconditioners in the implicit parallel accurate
  reservoir simulator (IPARS)}.
\newblock {\em Numerical Linear Algebra with Applications}, 8(8):537--549, dec
  2001.

\bibitem{LIVESCU2010138}
S.~Livescu, L.~Durlofsky, K.~Aziz, and J.~Ginestra.
\newblock A fully-coupled thermal multiphase wellbore flow model for use in
  reservoir simulation.
\newblock {\em Journal of Petroleum Science and Engineering}, 71(3):138 -- 146,
  2010.
\newblock Fourth International Symposium on Hydrocarbons and Chemistry.

\bibitem{MJE05}
V.~Martin, J.~Jaffr{\'{e}}, and J.~E. Roberts.
\newblock {Modeling fractures and barriers as interfaces for flow in porous
  media}.
\newblock {\em SIAM Journal on Scientific Computing}, 26(5):1667--1691, 2005.

\bibitem{MMB2007}
S.~K. Matthai, A.~A. Mezentsev, and M.~Belayneh.
\newblock {Finite element - node-centered finite-volume two-phase-flow
  experiments with fractured rock represented by unstructured hybrid-element
  meshes}.
\newblock {\em SPE Reservoir Evaluation {\&} Engineering}, 10(06):740--756, dec
  2007.

\bibitem{MF07}
J.~E. Monteagudo and A.~Firoozabadi.
\newblock {Control-volume model for simulation of water injection in fractured
  media: incorporating matrix heterogeneity and reservoir wettability effects}.
\newblock {\em SPE Journal}, 12(03):355--366, sep 2007.

\bibitem{NBFK2019}
J.~Nordbotten, W.~Boon, A.~Fumagalli, and E.~Keilegavlen.
\newblock Unified approach to discretization of flow in fractured porous media.
\newblock {\em Computational Geosciences}, 23:225--237, 2019.

\bibitem{Peaceman78}
D.~Peaceman.
\newblock {Interpretation of Well-Block Pressures in Numerical}.
\newblock {\em Reservoir Simulation Symposium Journal SEPJ}, pages 183--194,
  1978.

\bibitem{Peaceman83}
D.~Peaceman.
\newblock {Interpretation of Well-Block Pressures in Numerical Reservoir
  Simulation with Nonsquare Grid Blocks and Anisotropic Permeability}.
\newblock {\em Reservoir Simulation Symposium Journal SEPJ}, pages 531--543,
  1983.

\bibitem{RJBH06}
V.~Reichenberger, H.~Jakobs, P.~Bastian, and R.~Helmig.
\newblock {A mixed-dimensional finite volume method for two-phase flow in
  fractured porous media}.
\newblock {\em Advances in Water Resources}, 29(7):1020--1036, jul 2006.

\bibitem{SBN12}
T.~Sandve, I.~Berre, and J.~Nordbotten.
\newblock {An efficient multi-point flux approximation method for Discrete
  Fracture-Matrix simulations}.
\newblock {\em Journal of Computational Physics}, 231(9):3784--3800, may 2012.

\bibitem{SMW:2003}
R.~Scheichl, R.~Masson, and J.~Wendebourg.
\newblock {Decoupling and block preconditioning for sedimentary basin
  simulations}.
\newblock {\em Computational Geosciences}, 7(4):295--318, 2003.

\bibitem{refthermo}
E.~Schmidt.
\newblock {\em {Properties of water and steam in S.I. units}}.
\newblock Springer-Verlag, 1969.

\bibitem{shi2005}
H.~Shi, J.~A. Holmes, L.~J. Durlofsky, K.~Aziz, L.~Diaz, B.~Alkaya, and
  G.~Oddie.
\newblock Drift-flux modeling of two-phase flow in wellbores.
\newblock {\em SPE Journal}, 10(01):24--33, 2005.

\bibitem{TFGCH12}
X.~Tunc, I.~Faille, T.~Gallou{\"{e}}t, M.~C. Cacas, and P.~Hav{\'{e}}.
\newblock {A model for conductive faults with non-matching grids}.
\newblock {\em Computational Geosciences}, 16(2):277--296, mar 2012.

\bibitem{Wolf03}
C.~Wolfsteiner, L.~J. Durlofsky, and K.~Aziz.
\newblock Calculation of well index for nonconventional wells on arbitrary
  grids.
\newblock {\em Computational Geosciences}, 7(1):61--82, 2003.

\bibitem{tracer2016}
F.~Xing, R.~Masson, and S.~Lopez.
\newblock {Parallel Vertex Approximate Gradient discretization of
  hybrid-dimensional Darcy flow and transport in discrete fracture networks}.
\newblock {\em Computational Geosciences}, 2016.

\bibitem{Xing.ea:2017}
F.~Xing, R.~Masson, and S.~Lopez.
\newblock Parallel numerical modeling of hybrid-dimensional compositional
  non-isothermal darcy flows in fractured porous media.
\newblock {\em Journal of Computational Physics}, 345:637--664, sep 2017.

\end{thebibliography}
\end{document}